\documentclass[12pt,fleqn,emtex]{article}
\usepackage{times}
\usepackage{graphicx}

\usepackage[german, american]{babel}  %Deutsch

\usepackage{times}
\usepackage{graphicx}
\usepackage{amsfonts}       %
\usepackage{amsmath} %[leqno]=Gleichungen werden links numeriert;
\usepackage{amssymb}        %
\usepackage{amstext}        %
\usepackage{amsthm}         %theorem-umgebungen
\usepackage[matrix,curve]{xypic}

\usepackage{epsfig}
\psfigdriver{dvips}
\usepackage{latexsym}
\usepackage[matrix,curve]{xypic}
\usepackage{rotating}
\rotdriver{dvips}

\newcommand{\Z}{\mathbb{Z}}

\newtheorem{Lemma}{Lemma}
\newtheorem{Theorem}{Theorem}
\newtheorem{Corollary}{Corollary}
\newtheorem{Proposition}{Proposition}

\DeclareMathSizes{12}{11}{8}{5} 
\usepackage{xypic,dsfont}

\begin{document}
%\section{}
%\section{}

\centerline{\Large\bf Model structures, categorial quotients
and } \centerline{\Large\bf representations of super commutative Hopf
algebras II} \centerline{\Large\bf The case $\bf Gl(m\vert n)$}

%\bigskip\noindent
%\centerline{Preliminary version, 7.10.2010}

\bigskip\noindent
\centerline{R.Weissauer}

\goodbreak
\bigskip\noindent
%\section{\bf Weights and representations}

\bigskip\noindent

\section{Introduction} 

\bigskip\noindent
Let $k^{m\vert n}$ be a $\Z/2$-graded vector space $X=k^m \oplus k^n$ over a field $k$ of characteristic zero with even part $k^m$ and odd part $k^n$. The super linear group $G=Gl(m\vert n)$ contains the classical reductive algebraic $k$-group $Gl(m)\times Gl(n)$. The Lie super algebra
of $G$ is $Lie(G)=End_k(k^m\oplus k^n)$, which decomposes into the subspace
of endomorphisms which preserve respectively do not preserve the $\Z/2\Z$-grading.
The 
even part of $Lie(G)$ can be identified with $Lie(Gl(m)\times Gl(n))$. 
The Lie super bracket on $End_k(k^m\oplus k^n)$ is defined by $[X,Y] = X\circ Y - (-1)^{\vert X\vert \vert Y\vert} Y\circ X$ for graded endomorphisms $X,Y$ in $End_k(k^m\oplus k^n)$. We suppose $m\geq n$.

\bigskip\noindent
An algebraic representation of $G$ over $k$ is a homomorphism $$\rho: Gl(m)\times Gl(n) \longrightarrow Gl(V)$$ of algebraic groups over $k$, where $V=V_+ \oplus V_-$ is a finite dimensional $\Z/2\Z$-graded $k$-vectorspace, 
together with a $k$-linear map
$$ Lie(\rho): Lie(G) \longrightarrow End(V) $$
so that 
\begin{enumerate}
\item The parity on $V$ is defined by the eigenspaces of $\rho(E_m,-E_n)$,
\item $Lie(\rho)$ is  parity preserving for the natural $\Z/2\Z$-grading on $End_k(V)$ induced by 
$V=V_+\oplus V_-$, 
\item $Lie(\rho)$ is $\rho$-equivariant and coincides with the Lie derivative of $\rho$ on the 
even part $Lie(Gl(m)\times Gl(n))$ of $Lie(G)$,
\item $Lie(\rho)$ respects the Lie super bracket.
\end{enumerate}
A one dimensional representation $\rho$ of $G$, which is the analog of the determinant, is the Berezin $Ber_{m\vert n}$
defined by $\rho(g_1\times g_2)= det(g_1)/det(g_2)$ so that $Lie(\rho)$ is the super trace on $End(k^{m\vert n})$. The representation space of the Berezin
is $k^{1\vert 0}$ or $k^{0\vert 1}$ depending on $n$ modulo two.

\bigskip\noindent
Let $\cal T$ denote the abelian $k$-linear tensor category of algebraic representations of $G$.  As a $k$-linear abelian category $\cal T$ decomposes into a direct sum of blocks $\Lambda$. 
We show that there exists a purely transcendent field extension $K/k$
of transcendency degree $n$ and a $K$-linear weakly exact tensor functor (see section \ref{semisim})
$$ \varphi:\  {\cal T}\otimes_k K  \ \longrightarrow  \ sRep_K(H) \quad , \quad H=Gl(m-n) $$
from the $K$-linear scalar extension of $\cal T$ to the semisimple category of finite dimensional algebraic super representations
of the reductive algebraic $K$-group $H=Gl(m-n)$ defined over $K$. The simple objects 
of $sRep_K(H)$ are  the $\rho[i]$ for the irreducible algebraic representations $\rho$ of $H$, up to a parity 
shift of the grading ($i=0$ or $i=1$). We prove in corollary \ref{Wakimoto} that the image of a simple object $V$ in $\cal T$
becomes zero under the functor $\varphi$ if and only if $V$ is a simple object, which is not maximal atypical (see section \ref{weights}). On the other hand, in the main theorem of section \ref{mainth} we prove that for maximal  atypical simple
objects the image $\varphi(V)$ is an isotypic representation
$m(V) \cdot \rho(V)[p(V)]$ in $sRep_K(H)$, where the multiplicity $m(V)$ is $>0$ and where
$\rho(V)$ is an irreducible representation of $H$ which only depends on the block $\Lambda$ of $V$. The parity shift $p(V)$ can be easily computed.
The computation of the multiplicity $m(V)$ is more subtle. We show $$1 \ \leq \ m(V)\ \leq \ n! \ .$$ If $V$ is simple with  highest weight $\mu$,
we prove in section
\ref{Weyl} the recursion formula
$$  m(\mu) \ = \  { n \choose n_1,\ldots,n_r } \cdot\ \prod_{i=1}^r \ m(\mu_i) $$
for $m(V)=m(\mu)$, which allows to express $m(\mu)$ in terms of multinomials coefficients and multiplicies
$m(\mu_i)$ for smaller $n$. Via 
$$ sdim_k(V) \ = \ (-1)^{p(V)} \cdot m(V) \cdot dim_k(\rho(V)) $$
this formula for $m(\mu)$ gives a rather explicit formula 
for the super dimension of a maximal atypical irreducible
representation $V$ with weight $\mu$, 
since the classical Weyl dimension formula for $Gl(m-n)$ 
computes $dim_k(\rho(V))$.

\bigskip\noindent

\section{Weights} \label{weights}

\bigskip\noindent
Let $k$ be a field of characteristic zero, and let $G= Gl(m\vert n)$ denote the superlinear
group over $k$. In the following we always assume $m\geq n$.

\bigskip\noindent
In this section we review some fundamental facts on highest weights. For more details see [BS1] and [BS4]. Let 
${\cal T}$ denote  the $k$-linear rigid 
tensor category $Rep_k(\mu,G)$ of $k$-linear algebraic super representations $\rho$ of $G$ on $\Z/2\Z$-graded finite dimensional
$k$-super vector spaces, such that $\rho(id_m,-id_n)$ induces the parity automorphism
of the underlying $\Z/2\Z$-grading of the representation space of $\rho$. Here 
$  \mu: {\Z}/2{\Z} \mapsto (id_m,-id_n) \in Gl(m)\times Gl(n) \subset Gl(m\vert n)$
as in [BS4].  The category $\cal T$ admits a $k$-linear anti-involution ${}^*\! : \cal T\to T$ so that $A^* \cong A$ holds for all simple objects $V$ and all simple projective objects $V$ of $\cal T$.
The intrinsic dimension $\chi(V)$
 in a rigid tensor category $\cal T$ with $End_{\cal T}({\mathbf 1}) \cong k$ is $\chi(A)= eval_V\circ coeval_V$. It is preserved by tensor functors. In our case $\chi(V)$ is the super dimension $sdim_k(A) = dim_k(V_+) - dim_k(V_-)$ of the underlying super vectorspace $V=V_+\oplus V_-$. This is easily seen using the forget functor ${\cal T}\to svec_k$  to the category $svec_k$ of finite dimensional $k$-super vector spaces.

\bigskip\noindent
The isomorphism classes $X^+$ of the irreducible finite dimensional representations of $Gl(m\vert n)$ are indexed by their highest weights $\lambda=(\lambda_1,..\lambda_m; \lambda_{m+1},.. ,\lambda_{m+n})$. Here  $\lambda_1 \geq ... \geq \lambda_m$ and $\lambda_{m+1} \geq ... \geq \lambda_{m+n}$ are integers, and every $\lambda\in {\Z}^{m+n}$ satisfying these inequalities
occurs as the highest weight of an irreducible representation $L(\lambda)$. 
The trivial representation $\mathbf 1$ corresponds to $\lambda =0$.
Assigned to each highest weight $\lambda\in X^+$  are two subsets of the numberline $\bf Z$, namely the set
$$ I_\times(\lambda)\ =\ \{ \lambda_1  , \lambda_2 - 1, .... , \lambda_m - m +1 \} $$
 of cardinality $m$, respectively the set of cardinality $n$
$$ I_\circ(\lambda)\ = \ \{ 1 - m - \lambda_{m+1}  , 2 - m - \lambda_{m+2} , .... , n-m - \lambda_{m+n}  \} \ .$$
Following the notations of [BS4]
the integers in $ I_\times(\lambda) \cap I_\circ(\lambda) $ are labeled by $\vee$, the remaining ones in $I_\times(\lambda)$ resp. $I_\circ(\lambda)$ are labeled by $\times$ resp. $\circ$. All other integers are then labeled by  $\wedge$. 
This labeling of the numberline $\Z$ uniquely characterizes the weight vector $\lambda$. If the label $\vee$ occurs $r$ times in the labeling, then $r$ is called the degree of
atypicality of $\lambda$. Notice that $0 \leq r \leq n$, and $\lambda$ is called
maximal atypical if $r=n$. Let $\cal A \subset T$ denote the full abelian subcategory
generated by the representations $L(\lambda)$ for maximal atypical weights $\lambda$.
In the following we usually identify $X^+$ with the set of all labelings of the numberline, such that
$\vee$ occurs $r$ times for some $0\leq r\leq n$ and $\times$ respectively $\circ$ occurs $m-r$ respectively $n-r$ times. There are two natural orderings on $X^+$, the Bruhat ordering and the coarser weight ordering.
For $\lambda\in X^+$ let ${\cal T}^{\leq \lambda}$ respectively ${\cal T}^{< \lambda}$
denote the full subcategories of $\cal T$ generated by objects, all whose Jordan-H\"older
constituents are simple modules $L(\mu)$ with highest weights $\mu\leq \lambda$
(resp. $\mu < \lambda$) with respect to the weight ordering.

%\bigskip\noindent
%{\it Example}. Assume $m\geq n$. The trivial representation gives
%$$ \vee\ \mbox{ at } \quad -m+1,...,n-m $$
%$$ x \ \mbox{ at } \quad n-m+1,...,0 $$
%$$ \circ\ \mbox{ nowhere }  $$
%$$ \wedge\ \mbox{ elsewhere } \ .$$
%Moving $\vee$ to the left, gives $\lambda=(0,0,-1,1,0)$.
%One further move to the left gives $\lambda=(0,0,-2,2,0)$
%and so on.

\bigskip\noindent
The abelian category $\cal T$ decomposes into blocks $\Lambda$, defined by the eigenvalues of a certain elements in the center of the universal enveloping algebra of the
Lie superalgebra $gl(m\vert n)$. Two irreducible representations $L(\lambda)$ and $L(\mu)$ are in the same block if and only if the weights $\lambda$ and $\mu$ define labelings with the same position of the labels $\times$ and $\circ$. The degree of atypicality is a block invariant, and
the blocks $\Lambda$ of atypicality $r$ are in 1-1
correspondence with pairs of disjoint subsets of
$\Z$ of cardinality $m-r$ resp. $n-r$. The irreducible
representations of each block $\Lambda$ are in 1-1 correspondence
to the subsets of cardinality $r$ in the numberline with the subset of all $\times$ and all $\circ$
removed.

\bigskip\noindent
{\it Example}. The Berezin representation $Ber=Ber_{m\vert n}$ of $Gl(m\vert n)$ has highest weight
$\lambda=(1,..,1,1..,1;-1 ,..., -1)$ with $m$ digits $1$ and $n$
digits $-1$. Its superdimension is $sdim_k(Ber_{m\vert n}) = (-1)^n$ and
its dimension is $1$. All powers $Ber^k$ for $k\in\Z$
of the Berezin are maximal atypical, and $L(\mu)\in \cal A$ iff
$L(\mu) \otimes Ber^k \in \cal A$.

\bigskip\noindent
\section{Maximal atypical blocks $\Lambda$}

\bigskip\noindent
A block $\Lambda$ is maximal atypical if and only if it does not contain any label $\circ$. 
Assume $\Lambda$ is maximal atypical. Let $j$ then be minimum of the subset of all $x$ defined by $\Lambda$, or $j=1$ if there is no cross.
The uniquely defined weight $\lambda$, where the labels $\vee$ are at the positions
$j-1,...,j-n$ is called the {\it ground state} of the block. There are higher ground states
for $N=0,1,2,3,..,$ where the  labels $\vee$ are at the positions
$j-1-N,...,j-N-n$. The corresponding weight vectors of these ground states
are
$$ \lambda_N\ =\ (\lambda_1,...,\lambda_{m-n},\lambda_{m-n}-N,...,\lambda_{m-n}-N; -\lambda_{m-n}+N,...,-\lambda_{m-n}+N) \ $$
where $\{ \lambda_1,\lambda_2 -1,...,\lambda_{m-n}+1-m+n \}$ gives the positions of the
labels $\times$. For $m=n$ these are the powers $Ber^{-N}$ of the Berezin, and the ground state is the trivial representation $\mathbf 1$.  All ground states define irreducible representations $L(\lambda_N)$ in the given maximal
atypical block $\Lambda$.  For $\mu_i=\lambda_i - \lambda_{m-n}$, the representation
$$  L(\lambda_0) \otimes Ber^{-\lambda_{m-n}} \ =\ L(\mu_1,...,\mu_{m-n},0,...,0;0,...,0) \ ,$$
 again is an irreducible maximal atypical representation
(usually in another block). It is covariant in the following sense:

\bigskip\noindent
{\it Covariant representations}.  Let 
$\lambda_1 \geq \lambda_2 \geq ... $ be any partition $\lambda$
of some natural number $N=deg(\lambda):=\sum_{\nu} \lambda_\nu$.
Associated to this partion is the
covariant representation $$\{ \lambda \}\ := \ Schur_{\lambda}(k^{m\vert n})\ $$
as a direct summand defined by a Schur projector of the $N$-fold tensor product $X^{\otimes N}$ of the
standard representation $X$ of $Gl(m\vert n)$ on $k^{m\vert n}$.
The representation $\{ \lambda \}$ so defined is zero iff $\lambda_{m+1} > n$, and is nontrivial and irreducible
otherwise ([BR]). If the {\it hook condition} $\lambda_{m+1}\leq n$ is satisfied,
we may visualize this by considering the Young diagram attached to $\lambda$ with
first column $\lambda_1$, second column $\lambda_2$ and so on.
Let $\beta$ denote the intersection of the Young diagram with the box $\{ (x,y) \vert x \leq m, y \leq n \}$. Then $\lambda$ has the following shape
with subpartitions $\alpha,\beta,\gamma$ obtained by intersection we the three hook sectors 
$$ \lambda\ =\ \begin{matrix} \alpha &   \cr \beta & \gamma \cr  & \cr\end{matrix} \ .$$ 
The transposed $\gamma^*$ of the partition $\gamma$ is again
a partition with $(\gamma^*)_i=0$ for $i>n$. We quote from [BR] and [JHKTM] the assertion of the next 

\bigskip\noindent
\begin{Lemma} \label{cov}{\it If $\{ \lambda \}$ is not zero, then $\{ \lambda
\} \cong L(\mu)$ is irreducible with highest weight $\mu$ defined
by $\mu_i = \lambda_i$ for $i=1,...,m$, and $\mu_{m+i} =
max(0,(\lambda^*)_i-m)=(\gamma^*)_i$ for $i=1,...,n$. In other words
$$ \Biggl\{ \begin{matrix} \alpha &   \cr \beta & \gamma \cr   \end{matrix}\Biggr\}
\ \cong \ L\biggl( \begin{matrix}  \alpha   \cr \beta \cr
\end{matrix} \  ;  \begin{matrix}      \cr\gamma^* \cr
\end{matrix} \biggr) \ .$$} \end{Lemma}

\bigskip\noindent
This implies

\begin{Lemma} \label{atyp} A covariant representation $\{ \lambda\}$ 
attached to a partion $\lambda$ with the hook condition is maximal atypical if and only if
$\lambda_{m-n+1}=0$, and then $\mu=\lambda$ holds and
$$ \{ \lambda\}\ =\ L(\lambda_1,...,\lambda_{m-n},0,..,0;0,..,0) \ .$$
\end{Lemma}

\bigskip\noindent
{\it Proof}. One direction is clear. For $\lambda_{m-n+1}=0$ the representation $\{\lambda\}$ corresponds to the ground state of some maximal atypical block as explained above.
For the converse assertion notice that $I_\times(\mu) \subset [1-m,..,\infty)$ for the highest weight
$\mu$ of  the representation $\{\lambda\}$ by lemma \ref{cov}. Hence  $1-m-\mu_{m+1}$
is in $I_\circ(\mu)$, but not in $I_\times(\mu)$ if $\mu_{m+1} = (\gamma^*)_1 >0$. Hence, if 
$L(\mu)$ is maximal atypical, we conclude $\gamma^*=0$ and hence $\gamma=0$.
So we can assume $m>n$. Then $I_\circ(\mu) = [1-m,2-m,...,n-m]$, since $\gamma^*=0$. Therefore none of the $\mu_1, \mu_2-1,..,\mu_{m-n}+n-m+1$ is in $I_\circ(\mu)$, since $\mu_i = \lambda_i \geq 0$ for $i=1,..,m$. If $L(\mu)$ is maximal atypical, then all remaining $\mu_i + 1- i $ for $i=m-n+1,..,m$ must be contained in $I_\circ(\mu)$.
Since $\lambda_{m-n+1}=\mu_{m-n+1}$ is $\geq 0$, this implies $\lambda_{m-n+1} + n- m \in I_0(\mu)$
and therefore  $\lambda_{m+n+1}=0$. QED

\bigskip\noindent
Consider 
$$ \Pi\ =\ \Lambda^{m-n}(X) \otimes Ber_{m\vert n}^{-1}\ =\ L(0,..,0,-1,..,-1;1,..,1) $$
with $n$ digits $1$ and $-1$, the first higher ground state in the $\mathbf 1$-block. 

\begin{Lemma} \label{ground} Let $\Lambda$ be a maximal atypical block. For  the
ground states $L(\lambda_N)$ of order $N=0,1,2,..$ in this block $\Lambda$ 
we obtain 
$$  L(\lambda_{N}) \otimes \Pi \ \cong \ L(\lambda_{N+1}) \oplus R_N $$
for certain $R_N$, whose projection to all maximal atypical blocks of $\cal T$ are zero. 
\end{Lemma}

\bigskip\noindent
{\it Proof}. By a twist with $Ber^{\lambda_{m-n}-N-1}$ we can easily reduce to the case $N=0$ and
$\Lambda_{m-n+1}=0$. Then $\lambda_N=\lambda$ for 
$ \lambda=(\lambda_1,...,\lambda_{m-n},0,..,0; 0,..,0) $
defines a covariant representation $L(\lambda_N)=\{\lambda\}$. By the well known properties of Schur projectors, 
$$ L(\lambda_N) \otimes \Lambda^{m-n}(X) \ = \ \bigoplus_{\rho} \ [\rho: \lambda,\Lambda^{m-n}] \cdot 
\{ \rho \}$$ holds with the Littlewood-Richardson coefficients $[\rho: \lambda,\Lambda^{m-n}]$.
It is well known that $[\rho: \lambda,\Lambda^{m-n}]\neq 0$ implies $\rho_{m-n+1}>0$
unless $\rho=\lambda+\Lambda^{m-n}$. Hence by lemma \ref{atyp}
all summands $\{\rho \}$ are not maximal atypical, except for $\{\rho\} = \{\lambda + \Lambda^{m-n}\}$. By twisting with $Ber^{-1}$ our claim follows. QED

\bigskip\noindent

\bigskip\noindent
\section{\bf The stable category $\cal K$}

\bigskip\noindent
The abelian category ${\cal T=\cal T}_{m\vert n}$  is a Frobenius category, i.e. it has enough projective objects and the injective and projective objects coincide. Let ${\cal K=\cal K}_{m\vert n}$ be the quotient category
with the same objects as $\cal T$, but with $Hom_{{\cal K}}(A,B)$ defined as the 
quotient of $Hom_{{\cal T}}(A,B)$ by the $k$-subvectorspace of all homomorphisms which factorize over a projective module. The natural functor $$\alpha: {\cal T}\longrightarrow {\cal K}$$ is a $k$-linear tensor functor. The category ${\cal K}$ is a triangulated category with a suspension functor $S(A)=A[1]$, and the quotient functor $\alpha$ associates to exact sequences in $\cal T$ distinguished triangles in ${\cal K}$. Furthermore
$$ Ext_{\cal T}^i(A,B)\ \cong\ Hom_{{\cal K}}(A,B[i]) \quad , \quad \forall \ i> 0 \ .$$ 
For simple  $X$ in $\cal T$ either $X$ is projective in $\cal T$ and hence
zero in ${\cal K}$, or $X$ is not zero in ${\cal K}$ and
$$ Hom_{{\cal K}}(X,X) \ = \ k \cdot id_X \ $$
is one dimensional.
Notice that $A[1] \cong I/A$ for an embedding $A \hookrightarrow
I$ into a projective object $I$, and similarly $A[-1] \cong Kern(P \to A)$ for a projective resolution $P\to A$. Since $P^* \cong P$ and $\cal T$ is a Frobenius category, this implies that ${}^*$ induces an involution
of the stable category ${\cal K}$ such that $(A[n])^* \cong A^*[-n]$.

\bigskip\noindent
\begin{Theorem} \label{ext} ([BS2] corollary 5.15). $dim_k(Ext_{\cal T}^i(L(\lambda), L(\mu)))$ is equal to
$$\sum_{j+k=i}\ \sum_{\nu} \ \ dim_k(Ext_{\cal T}^j(V(\nu), L(\lambda))) \cdot dim_k(Ext_{\cal T}^k(V(\nu), L(\mu))) \ .$$ 
\end{Theorem}

\bigskip\noindent
\begin{Theorem} ([BS2], corollary 5.5).
$dim_k(Ext_{\cal T}^i(V(\lambda), L(\mu)))=0$ unless 
$\lambda\leq \mu$ and 
$i \equiv l(\lambda,\mu)$
modulo 2. 
\end{Theorem}

\bigskip\noindent
Here $l(\lambda,\mu)$ denoted the minimum number of transpositions of neighbouring $\vee\!\wedge$ pairs needed to get from $\lambda$ to $\mu$, where neighbouring means separated only by $\circ$'s and $x$'s. Put $p(\lambda)=\sum_{i=1}^n \lambda_{m+i}$.
If $\lambda$ and $\mu$ are maximal atypical, then
$p(\lambda)\equiv p(\mu) + l(\mu,\nu)$ modulo two. Indeed it suffices to show this in the case of
neighbours $\lambda$ and $\mu$ where $l(\lambda,\mu)=1$. For
maximal atypical weights $I_x \cap I_0 = I_0$,  a single transposition modifies $\sum_{j \in I_x\cap I_0} j = -\sum_{i=1}^{n} \lambda_{m+i} + \sum_{i=1}^{n} (i-m) $ by one. 
Hence the last two theorems imply

\bigskip\noindent
\begin{Lemma} \label{parity} For $L(\lambda)$ and $L(\mu)$ in $\cal A$ and $i \geq 0$ we have 
$Hom_{{\cal K}}(L(\lambda), L(\mu)[i])=0$ unless $p(\lambda) \equiv p(\mu) + i$
modulo two.
\end{Lemma}

\begin{Lemma} \label{trivial}  Let $\cal A$ be a $k$-linear category
and let $A$ and $B$ be objects of $\cal A$ such that
$End_{\cal A}(A)\cong k$ and $End_{\cal A}(B) \cong k$.
Let $\varphi: B\to A$ and $\psi:A\to B$ be morphisms in $\cal A$.
Then, either $\varphi$ and $\psi$ are isomorphisms, or
$\varphi\circ\psi=0$ and $\psi\circ\varphi =0$.
\end{Lemma}

\bigskip\noindent
{\it Proof}. Suppose $\varphi\circ\psi\neq 0$. Then by a rescaling
we can assume $\varphi \circ \psi = id_A$.  Then $\varphi \circ \psi \circ \varphi = \varphi$.
Since $\psi \circ \varphi = \lambda \cdot id_B$
for some $\lambda\in k$ and since $\varphi\neq 0$ by our assumption, we conclude $\lambda =1$.
Hence $\varphi$ and $\psi$ are isomorphisms.
The same conclusions hold if
 $\psi\circ\varphi \neq  0$. QED

\bigskip\noindent
\section{Kostant weights}\label{Kosw}

\bigskip\noindent
By [BS2], lemma 7.2 a weight $\mu$ is a Kostant weight, i.e. 
satisfies for every $\nu\in \Lambda$
$$ \sum_{i=0}^\infty \ dim_k( Ext^i_{\cal T}(V(\nu),L(\mu)) \ \leq \ 1 \ , $$
if and only if no subsequence of type $\vee\! \wedge\! \vee\! \wedge$ occurs in its labeling.
All ground states $\lambda_N$ of the maximal atypical blocks are Kostant weights.
In particular there is at most one index $i\! =\! i(\nu)$, depending on $\lambda_N$ and $\nu$, such that $Ext^i_{\cal T}(V(\nu),L(\lambda_N)) \neq 0$, and in this case $Ext^{i(\nu)}_{\cal T}(V(\nu),L(\lambda_N)) \leq 1$.

\bigskip\noindent
By [BS2] corollary 5.5 and the complementary formula (5.3) for $p_{\nu,\mu}$
in loc. cit. $
Ext^i_{\cal T}(V(\nu),L(\mu))=0$ holds unless $\nu,\mu$ are in the same block $\Lambda$, unless $\nu \leq \mu$ in the Bruhat ordering
and $i \leq l(\nu,\mu)$. Suppose $\nu,\mu$ are in the same block $\Lambda$, and
suppose $\nu\leq \mu$ holds in  the Bruhat ordering.  Then for  $i=l(\nu,\mu)$
we have $Ext_{\cal T}^i(V(\nu),L(\mu))\neq 0$ by inspection of the  formula (5.3) in [BS2].
Indeed in the set $D(\nu,\mu)$ of labeled cap diagrams $C$
defined in loc. cit., there exists at least one cap diagram with $\vert C\vert =0$, 
since $\nu\leq \mu$ in the Bruhat ordering implies $l_i(\nu,\mu)\geq 0$ for all $i\in I(\Lambda)$
in the notations of loc. cit. Since the leftmost vertex of a small cap is always $\vee$ and hence is contained in $I(\Lambda)$, there exists some $C\in D(\nu,\mu)$ with $\vert C\vert=0$.
This implies

\begin{Lemma} \label{Kostant}
For $L=L(\mu)$ and a Kostant weight $\mu$ the $k$-vectorspace $Ext^i_{\cal T}(V(\nu),L)$
is one dimensional, if $\nu,\mu$ are in the same block $\Lambda$ and $\nu\leq \mu$ holds in the Bruhat ordering and $i=l(\nu,\mu)$, and it is zero otherwise. 
\end{Lemma}

\bigskip\noindent
If $\mu=\lambda_N$ is one of the ground state weights of the block $\Lambda$,
then the conditions $\nu\in \Lambda$, $\nu\leq \mu$ and $i=i(\nu,\mu)$ only
depend on the relative positions of the labels $\vee$ in the numberline after the crosses
defined by the block $\Lambda$ are removed. So they do not depend on the block, but
only on the number $n$ of labels $\vee$ of the Kostant weight. 
This implies

\begin{Corollary} \label{dimen}
For any block $\Lambda$ of $\cal T$ and ground state representation
$L=L(\lambda_N)$ of this block we have for all $j\geq 0$ 
$$ dim_k(Ext^j_{\cal T}(L,L))\ =\ dim_k(Ext^j_{\cal T}({\mathbf 1},{\mathbf 1})) \ .$$
\end{Corollary}

\bigskip\noindent
{\it Proof}. 
By theorem \ref{ext} we get $dim_k(Ext_{\cal T}^j(L,L) =
\sum_\nu dim_k(Ext_{\cal T}^{i(\nu)}(V(\nu),L))^2 $ with summation 
over all $\nu$ such that $j=2i(\nu)$. Since by lemma \ref{Kostant} the summation conditions
and the dimensions in this sum only depend
on $j$ and not on the chosen ground state $L$ or block $\Lambda$, 
we may replace $L$ by the ground state ${\mathbf 1}$ of the trivial block. QED

\bigskip\noindent 

\bigskip\noindent
\section{\bf The localization $\cal B$ of the tensor category}

\bigskip\noindent
For the moment let $\cal K$ be any $k$-linear triangulated tensor category (meaning symmetric monoidal).
Then there exists a triangulated tensor functor ${\cal K \to \cal K}^\sharp$ of idemcompletion
(see [BS]).

\bigskip\noindent
Let $k_{\cal K}=End_{\cal K}({\mathbf 1})$ be the central $k$-algebra of $\cal K$.
Let $u\in \cal K$ be an invertible element (we use this for $u={\mathbf 1}[1]$).
The monoidal symmetry $\sigma_u: u\otimes u \cong u\otimes u$ is given by multiplication
with an element $\epsilon_u \in k_{\cal K}^*$ with $(\epsilon_u)^2=1$.
Furthermore $$R^\bullet_{\cal K} \ = \ \bigoplus_{i\in\Z} \ Hom_{\cal K}({\mathbf 1},u^{\otimes i} )$$
becomes a supercommutative ring with the parity $\epsilon_u$,
i.e. $f g = (\epsilon_u)^{ij} g f$ for homogenous elements of degree $i$ and $j$.
See [Ba], Prop. 3.3.
Let $R \subset R^\bullet_{\cal K}$ be the graded subring generated by the
elements of degree $\geq 0$.
Let $S \subset R^\bullet_{\cal K}$ be a multiplicative (and even, if $\epsilon_u \neq 1$) subset,
then the ring localization $S^{-1}R^\bullet_{\cal K}$ is defined.
Define a new category by the degree zero elements
$$ Hom_{S^{-1}\cal K}(A,B) : = (S^{-1}  Hom_{\cal K}^\bullet(A,B))^0 $$
of the localization of the graded $R^\bullet_{\cal K}$ module
$$Hom_{\cal K}^\bullet(A,B) \ =\ \bigoplus_{i\in\Z} \ Hom_{\cal K}(A,u^{\otimes i} \otimes B)\ .$$
For $M\in \cal K$ the annulator of $Hom_{\cal K}^\bullet(M,M))$ in $R$
is a graded ideal, which defines the support variety $V(M)$ as the spectrum of the quotient ring.

\bigskip\noindent
\begin{Theorem}([Ba], thm. 3.6). {\it As a tensor category
$S^{-1}\cal K$ is equivalent to the Verdier quotient of the tensor
category $\cal K$ divided by the thick triangulated tensor ideal
generated by the cones of the morphisms in $S$. The quotient category
$\cal B$ is a $S^{-1}R^{\bullet}_{\cal K}$-linear category. The quotient functor is a 
$k$-linear triangulated tensor functor}. \end{Theorem}

\bigskip\noindent
In the following we only apply this for the stable category $\cal K$ of the representation category $\cal T$
for $u={\mathbf 1}[1]$. In this case $\epsilon_u=-1$.
Then we have

\bigskip\noindent
\begin{Proposition}(\mbox{[BKN1], 8.11}).\label{grad} {\it The graded ring $R$
of the Lie super algebra $psl(n\vert n)$ is isomorphic to a graded
polynomial ring $k[\zeta_2,\zeta_4,...,\zeta_{2n-2},\xi_n,\eta_n]$
of transcendency degree $n+1$. } \end{Proposition}

\bigskip\noindent
By using restriction from $Gl(m\vert n)$ for $m>n$ one finds that
for $Gl(n\vert n)$ the isomorphisms $\zeta_2,..,\zeta_{2n-2}$ also exist.
Now suppose $m=n$.
Then by proposition \ref{grad} there exists a power $\cal L$ of the Berezin such that
$\xi_n: {\mathbf 1} \to {\cal L}[n]$ and  $\eta_n: {\mathbf 1} \to {\cal L}^{-1}[n]$. The product
$\zeta_{2n}=\eta_n\xi_n$ is in $R$. In the appendix \ref{class} we show ${\cal L} = Ber_{n\vert n}$.

\bigskip\noindent
\begin{Proposition}([BKN2], p.23). \label{pol}{\it The graded ring $R$
of the category $Gl(m\vert n)$ is isomorphic to a graded
polynomial ring $S= k[\zeta_2,\zeta_4,...,\zeta_{2n}]$ of
transcendency degree $n$ for all $m\geq n$.} \end{Proposition}

\bigskip\noindent
For the support variety $V(M)$ of an object $M\in \cal K$, being defined as above,
we quote from [BKN1], section 7.2, p. 29 and [BKN2], 4.8.1

\bigskip\noindent
\begin{Theorem}([BKN2], thm. 4.8.1). \label{1}{\it For $Gl(m\vert
n)$ the dimension of the support variety of a simple object
$L(\lambda)$ is equal to the degree of atypicality of $L(\lambda)$.}
\end{Theorem}

\bigskip\noindent
\begin{Theorem}([BKN2], thm. 4.5.1 and 4.8.1). \label{2} {\it For $Gl(m\vert n)$
%\bigskip\noindent
%{\bf Theorem}([BKN1], corollary 6.4.1).{\it
the  support variety of a simple maximal atypical object
$L(\lambda)$ of $\cal T$ is $Spec(R)$.} \end{Theorem}

%\bigskip\noindent
%{\bf Theorem 3}([BKN2], cor. 3.3.1). {\it For $Gl(m\vert n)$ the support
%variety of a Kac module (or anti Kac module) is the point $0$.}

\bigskip\noindent
Now fix the supergroup $Gl(m\vert n)$.
Put  $K=Quot(R)=k(\zeta_{2},...,\zeta_{2n})$, and put $S=R \setminus \{ 0\}$.
Notice, then as required $S$ only contains even elements by proposition \ref{pol} above.
Furthermore, since $R$ is an integral domain, $S^{-1}R$ is isomorphic to the extension field $K$ of $k$.
Let ${\cal B} =R^{-1}{\cal K}$ be the corresponding localized category.
If $\cal B$ is not idemcomplete, we replace it
by its idemcompletion $\cal B^\sharp$ from now on.
${\cal B}$ is a $K$-linear category.
It is obvious that the functor ${}^*$ respects $R$, hence induces a 
corresponding $K$-equivariant functor of ${\cal B}$. The natural quotient functor
$$   \beta: {\cal K} \ \longrightarrow \ {\cal B} \ $$
is a $k$-linear triangulated tensor functor.

%\bigskip\noindent
%{\bf Question}. {\it How does $*$ and $\vee$ act on $K$? Are they $K$-linear or $K$-equivariant functors?}

%\bigskip\noindent
%Obvious is $\zeta_2^* =\pm \zeta_2$ and $\zeta_4^* = \pm \zeta +
%const. \zeta_2^2$ etc. Since $\eta_n^* = \xi_n$ (without restriction of generality)
%the last $\zeta_{2n} = \xi_n\eta_n$ is invariant under ${}^*$.

\bigskip\noindent
\section{The homotopy category ${\cal H}$}

\bigskip\noindent
In the last section we defined the quotient category $\cal B$ of the stable
category ${\cal K}$. We now define another Verdier quotient category
$\cal H$ of ${\cal K}$, which is obtained by dividing ${\cal K}$ by the thick triangulated $\otimes$-ideal of ${\cal K}$ generated by the anti Kac modules. Recall that
for each highest weight $\lambda \in X^+$ there exists a cell module $V(\lambda)$ (or Kac module) in the sense of [BS4] in the category $\cal T$. We define the anti Kac modules by $V(\lambda)^*$
via the antiinvolution ${}^*$. 

\bigskip\noindent
\begin{Lemma}{\it For (anti)-Kac modules $V$ and $\zeta_{2i}\in R$
the space $\bigoplus_{n\in \Z} Hom_{\cal K}(V,V[n])$ is annihilated by a sufficiently high powers of $\zeta_{2i}$.}
\end{Lemma}

\bigskip\noindent
{\it Proof}. By [BKN2], thm. 3.2.1 the groups $Ext_{\cal T}^j(V,M)$ vanish for
fixed $M\in\cal T$ if $j>>0$.  Hence a high enough power of $\zeta_{2i}$ annihilates
$Hom_{\cal K}(V,V[n])$. By applying the functor ${}^*$ this carries over to $V^*$.
QED

\bigskip\noindent
Since $\zeta_{2i}$ becomes invertible in $\cal B$, this implies 
$Hom_{\cal B}(V,V)=0$ for each Kac module $V$. Hence

\bigskip\noindent
\begin{Lemma}\label{cell} {\it  The image
of a (anti)-Kac module $V$ is zero in $\cal B$.}
\end{Lemma}

\bigskip\noindent
Recall that in [W] we defined the homotopy category ${\cal H}$, which
is equivalent as a $k$-linear tensor category
to the Verdier quotient of the category ${\cal K}$ devided by the thick tensor ideal generated
by all anti-Kac modules. By the last lemma we obtain

\bigskip\noindent
\begin{Corollary} {\it The quotient functor $\beta: {\cal K} \to {\cal B}$ factorizes over the
homotopy quotient functor $\gamma: {\cal K} \to {\cal H}$.}
\end{Corollary}

\bigskip\noindent
We quote from [W] the following

\begin{Theorem} \label{hot} In the category $\cal H$  for simple objects
$M$ and $N$ in $\cal T$ with highest weights $\mu$ and $\lambda$ the following holds:
\begin{enumerate}
\item $End_{\cal H}(M) = k \cdot id_M$,
\item $Hom_{\cal H}(M,N)$ is a finite $k$-vectorspace,
\item $Hom_{\cal H}(V,N) \cong Hom_{\cal T}(V,N)$ for every cell module
$V=V(\mu)$. 
\item $Hom_{\cal H}(M,N) =0$ if $\mu < \lambda$ holds with respect to the weight ordering,
\item Let ${\cal H}^{\leq \lambda}$ denote the full subcategory quasi equivalent to the
image of ${\cal T}^{\leq \lambda}$ and similar ${\cal H}^{< \lambda}$ for ${\cal T}^{< \lambda}$ . Then suspension induces a functor
$$ [1]: \ {\cal H}^{\leq \lambda}\ \longrightarrow \ {\cal H}^{< \lambda} \ .$$
\end{enumerate}
\end{Theorem}

\bigskip\noindent
\begin{Lemma} \label{unit} $End_{\cal B}(1) = K$. 
\end{Lemma}

\bigskip\noindent
{\it Proof}. By [Ba], prop 3.3 we have $End_{\cal B}(1) = S^{-1}R^\bullet_{\cal K}$. Hence it suffices to show for $n>0$
that all morphisms
$$ \psi: {\mathbf 1} \longrightarrow {\mathbf 1}[-n] $$
are annihilated by a power of the element $\zeta_{2} \in R$.
Suppose $n=2i-1$ is odd. Then this is obvious, since $(\zeta_2)^i \circ \psi: {\mathbf 1}\to u={\mathbf 1}[1]$
vanishes by parity reasons. Indeed $Hom_{\overline {\cal T}}({\mathbf 1},{\mathbf 1}[1])=0$ by 
lemma \ref{parity}. Now suppose
$n=2i$ is even and consider $(\zeta_2)^i\cdot \psi$ in $End_{\cal K}({\mathbf 1})=k$. By the trivial lemma
\ref{trivial}
the morphism $\psi$ is an isomorphism unless $(\zeta_2)^i \cdot \psi =0$.
If $\psi$ were an isomorphism in $\cal K$, then also in $\cal B$ and therefore
also in $\cal H$. However $Hom_{\cal H}({\mathbf 1},{\mathbf 1}[-n])=0$ holds for all $n>0$ by the assertions 4) and 5) of theorem
\ref{hot}. QED

\bigskip\noindent
The same argument shows

\bigskip\noindent
\begin{Lemma}\label{nonvan}{\it For simple objects $M$ in $\cal T$
the image of $Hom_{\cal K}(M,M[-n])$  under the natural map
$$ Hom_{\cal K}(M,M[-n]) \longrightarrow  S^{-1} Hom^\bullet_{\cal K}(M,M)^0 = Hom_{\cal B}(M,M) $$ 
is zero for $n>0$.}
\end{Lemma}

\bigskip\noindent
{\it Proof}. $Hom_{\cal H}(M,M[-n])$ vanishes for simple $M$ and all $n>0$
by weight reasons. Use  part 4) and 5) of theorem \ref{hot}. QED

\bigskip\noindent
By lemma \ref{nonvan} for simple objects $M$ the endomorphism ring  $Hom_{\cal B}(M,M)$ is obtained by
quotients $f/r$ for $f:M\to M[i]$ for $f\in Hom_{\cal K}(M,M[i])$ and $r:{\mathbf 1}\to {\mathbf 1}[i]$ 
in $r\in Hom_{\cal K}({\mathbf 1},{\mathbf 1}[i])$ only for positive degrees $i\geq 0$,
of course modulo the usual equivalence defined by the localization $S^{-1}$.
We may suppose that the simple object $M$ is not projective, since otherwise
$M$ vanishes in $\cal K$ and hence in $\cal B$. Then $M\neq 0$ in $\cal M$, and
$Hom_{\cal T}(M,M) \cong Hom_{\cal K}(M,M) \cong k$. Hence
$Ann_R(Ext^\bullet_{\cal T}(M,M)) = Ann_R(\bigoplus_{i=0}^\infty Hom_{\cal K}(M,M[i]))$. 
 Hence
$r\in Ann_R(Ext^\bullet_{\cal T}(M,M))$ iff $r \cdot f =0$ for all $f\in 
Ext^\bullet_{\cal T}(M,M)$. This is related to
the support variety $V(M)$ of the simple object $M$ $$V(M) = Spec(R/Ann_R(Ext^\bullet_{\cal T}(M,M))) \ $$
as follows. There exists an $r\neq 0$ in $R$ that annihilates $\bigoplus_{i=0}^\infty Hom_{\cal K}(M,M[i])$ iff the support of $V(M)$ is not equal to $Spec(R)$.
The first statement is equivalent to $Hom_{\cal B}(M,M)=0$. Hence

\begin{Corollary}
Let $M$ be a simple object of $\cal T$. Then $M$ vanishes in $\cal B$ if and only if 
the support variety $V(M)$ is a proper subset of $Spec(R)$. 
\end{Corollary}

\bigskip\noindent
\begin{Corollary} \label{Wakimoto} {\it A simple object $M$ of $\cal T$ vanishes in $\cal B$
iff $M$ is not maximal atypical.} \end{Corollary}

\bigskip\noindent
{\it Proof}.  Use theorem \ref{1} and
\ref{2}.

\bigskip\noindent
\begin{Lemma}{\it  For simple objects $M$ and $N$ in $\cal T$
the space $Hom_{\cal B}(M,N)$ vanishes unless $M$ and $N$ have the same  parity in the sense of lemma \ref{parity}.} \end{Lemma}

\bigskip\noindent
{\it Proof}. As explained above, any morphism in $\cal B$ between $M$ and $N$ is of the form
$f/r$ for $f:M\to N[i]$ in $\cal K$, $r\in R$ for some even $i\geq 0$.
Hence $f$ corresponds to an element in $Ext_{\cal T}^i(M,N)$. We can
assume that $M$ and $N$ are maximal atypical, since otherwise $M$ and $N$ are zero in $\cal B$.  Then $Ext_{\cal T}^i(M,N)$ vanishes  by lemma \ref{parity} unless $M$ and $N$ have the same
parity. QED

\bigskip\noindent
We remark that $\zeta_{2}$ becomes an isomorphism in $\cal B$. Hence

\bigskip\noindent
\begin{Lemma} ${\mathbf 1}[2] \cong {\mathbf 1}$ in $\cal B$. \end{Lemma}

\bigskip\noindent

\section{The ground state categories ${\cal Z}_\Lambda$}\label{groundcat}

\bigskip\noindent
{\bf Definition}. Let $\Lambda$ be a maximal atypical block of $\cal T$ and $L$ be the ground state representation
of this block. Let
${\cal Z}_\Lambda$ denote the full subcategory of $\cal B$ of all objects isomorphic to a finite direct sum of $L$ and $L[1]$. For the unit block where $L={\mathbf 1}$ we simply write $\cal Z$.

\bigskip\noindent
\begin{Lemma} \label{thick} ${\cal Z}_\Lambda$
is a thick idemcomplete triangulated subcategory of $\cal B$,  i.e. it is closed under retracts and extensions
and the shift functor. 
\end{Lemma}

\bigskip\noindent
{\it Proof}. First suppose $L={\mathbf 1}$. Notice $$ Hom_{\cal B}(a\cdot {\mathbf 1} \oplus b\cdot {\mathbf 1}[1],
c\cdot {\mathbf 1} \oplus d\cdot {\mathbf 1}[1]) \ \cong \ Hom_{K}(K^a,K^c) \oplus Hom_K(K^b,K^d)$$ by lemma \ref{parity} and lemma \ref{unit}.
Hence by usual properties of matrix rings the category $\cal Z$ is an idempotent split category. By the same reasons $\cal Z$ is closed under retracts as well as under cones. Since ${\mathbf 1}[2]\cong {\mathbf 1}$, shifts preserve $\cal Z$. Hence $\cal Z$ is a thick idempotent split triangulated subcategory
of $\cal B$. The same carries over to ${\cal Z}_\Lambda$ by the next lemma.
QED

\begin{Lemma} Suppose $L$ is a ground state of a maximal atypical block $\Lambda$ 
in $\cal T$. Then $Ext^\bullet_{\cal T}(L,L)$ is a graded free module over $R=Ext^\bullet_{\cal T}({\mathbf 1},{\mathbf 1})$, hence in particular $$End_{\cal B}(L)\ \cong\ K \cdot id_L\ .$$ 
\end{Lemma}

\bigskip\noindent
{\it Proof}.  By theorem \ref{2} the annihilator of $Ext_{\cal T}^\bullet(L,L)$  in $R$ is trivial.
Hence the graded $R$-homomorphism $R\cdot id_L \to Ext_{\cal T}^\bullet(L,L)$ is injective.
By corollary \ref{dimen} the comparison of $k$-dimensions shows that
it is an isomorphism. QED

\bigskip\noindent
\begin{Lemma} \label{trivi} Suppose the image of an irreducible maximal atypical highest weight representation
$L(\lambda)$ in $\cal B$ is contained in $\cal Z$. Then $$ sdim_k(L(\lambda)) = (-1)^{p(\lambda)} \cdot m(\lambda) \ $$
for some integer $m(\lambda)>0$.
\end{Lemma}

\bigskip\noindent
{\it Proof}. By corollary \ref{Wakimoto} the object 
$L(\lambda)$ is not zero in $\cal B$. Lemma \ref{parity} together with the
assumption $
L(\lambda) \in \cal Z$ hence implies $L(\lambda) \cong c \cdot {\mathbf 1}[p'(\lambda)]$
for some $c\neq 0$ in $\Z$, where $p'(\lambda)\in\{0,1\}$ is uniquely defined.
Then $sdim_k(L(\lambda)) = (-1)^{p'(\lambda)} c$.
Two remarks. First $\chi_{\cal H}({\mathbf 1}[1]) = \chi_{\cal B}({\mathbf 1}[1])=-1$ for $u={\mathbf 1}[1]$, since $\epsilon_u = -1$. Notice, this holds in the homotopy category $\cal H$ by [W] and hence in $\cal B$.
Secondly for all $L(\lambda)$ in $\cal Z$ we have
$p'(\lambda) = p(\lambda) + const$, where $const \in Z$ is independent from $L(\lambda)$ in $\cal Z$ by lemma \ref{parity}. For $\lambda =0$ and $L(\lambda)=1$ this shows $const \in 2\Z$.  Therefore
$(-1)^{p'(\lambda)}=(-1)^{p(\lambda)}$.

\bigskip\noindent
Lemma \ref{ground} and corollary \ref{Wakimoto} imply
for the atypical representations $L(\lambda_N), L(\lambda_{N+1})$ and $\Pi$
the following relation in $\cal B$

\begin{Lemma} \label{N-state}
In $\cal B$ we have $L(\lambda_N) \otimes \Pi \cong L(\lambda_{N+1})$ for all $N\geq 0$.
\end{Lemma}

\section{The reductive group $Gl(m-n)$}\label{mainth}

\bigskip\noindent
Consider the tensor category ${\cal T}={\cal T}_{m\vert n}$ as before.
Let $H=Gl(m-n)$ denote the linear group over $k$ and let $Rep_k(H)$ denote the 
$k$-linear rigid semisimple tensor category of all algebraic
representations of $Gl(m-n)$ on finite dimensional $k$-vectorspaces.

\bigskip\noindent
Embedded in $Gl(m\vert n)$, with an immersion in the obvious way, is subgroup $H \times Gl(n\vert n)$. 
Restriction of a super representation of $Gl(m\vert n)$ on a finite dimensional
$k$-super vectorspace to the subgroup $H \times Gl(n\vert n)$ defines a $k$-linear
exact tensor functor $$Res: {\cal T}_{m\vert n} \longrightarrow Rep_k(H) \otimes_k
{\cal T}_{n\vert n}\ .$$ The restriction  of a projective
representation $P$ in ${\cal T}_{m\vert n}$ decomposes into a direct sum of isotypic
representations $P=\bigoplus_\rho P_\rho$ with respect to the action of  the reductive
group $H$. Each $P_\rho$ is a $Gl(n\vert n)$ representation, which is projective
as a direct summand of the projective object $P$ viewed as a super representation in ${\cal T}_{n\vert
n}$. Hence $Res({\cal P}) \subset Rep_k(H) \times {\cal
P}$. Similarly a morphisms in ${\cal T}_{m\vert n}$, which is
stably equivalent to zero, restricts to a direct sum of morphisms
$f_\rho$ with respect to the action of $H$,  such that each of the morphisms $f_\rho$
is stably equivalent to zero in ${\cal
T}_{n\vert n}$. Hence the restriction induces a tensor functor
$$ res: \ {\cal K}_{m\vert n} \ \longrightarrow \ Rep_k(Gl(m-n)) \otimes_k
{\cal K}_{n\vert n} \ .$$ The suspension $(.)[1]$ thereby maps to the
suspension $id_{Rep_k(H)} \otimes_k (.)[1]$, since an embedding
$X \hookrightarrow I $ decomposes into $Res(X)=\bigoplus_\rho
Res(X)_\rho \hookrightarrow \bigoplus_\rho Res(I)_\rho$. Thus
$res$ becomes a triangulated $k$-linear tensor
functor. The triangulated structure on  $Rep_k(Gl(m-n)) \otimes_k
{\cal K}_{n\vert n} $ is induced by the triangulated structure
on ${\cal K}_{n\vert n} $ in an obvious way, noticing that $Rep_k(H)$
is semisimple. 

\bigskip\noindent
Using detecting subalgebras as in [BKN] one can show that $Ext^n_{{\cal T}_{m\vert n}}(k,k)$ restricts properly and surjectively 
to $Ext^n_{{\cal T}_{n\vert n}}(k,k)$. By the universal property of the Verdier
quotient categories the
functor $res$ induces a functor from the Verdier quotient categories
${\cal B}_{m\vert n}$ of ${\cal T}_{m\vert n}$ 
$$ \gamma: {\cal B}_{m\vert n} \ \longrightarrow \ Rep_k(H) \otimes_k {\cal
B}_{n\vert n} \ .$$ 
This functor is a $K$-linear triangulated tensor functor.

\bigskip\noindent
Now use

\begin{Theorem} \label{Z}
For $m \geq n$ the image of the block of the trivial representation in the category ${\cal B}_{m\vert n}$ is equivalent as a $K$-linear triangulated tensor subcategory of ${\cal B}_{m\vert n}$
to the $K$-linear triangulated tensor category ${\cal Z} \sim svec_K$ of finite dimensional $K$-super vectorspaces.
\end{Theorem}

\bigskip\noindent
Taking this theorem for granted at the moment we proceed as follows:
We apply the last theorem for $m=n$, which allows us to consider $\gamma$ as a
functor $K$-linear triangulated tensor functor
$$ \gamma: {\cal B}_{m\vert n} \ \longrightarrow \ Rep_k(H) \otimes_k svec_K   \ .$$ 
The right side can be viewed as the category of finite dimensional $K$-algebraic
super representations of the reductive group $Gl(m-n)$ over $K$. Up to twist by powers of the determinant representation of $Gl(m-n)$ a basis of simple
objects is given by the representations $Schur_\mu(K^{m-n})$ and $Schur_\mu(K^{m-n})[1]$, where $\mu=\mu_1 \geq \mu_2 \geq ... \geq \mu_{m-n}\geq 0$ runs over the partions
of length $\leq m-n$.

\bigskip\noindent
Let us look what the functor $\gamma$ does with the image in ${\cal B}_{m\vert n}$ of the standard representation
$X_{m\vert n} = k^{m\vert n}$ of $Gl(m\vert n)$.
This standard representation
restricts to $Res(X_{m\vert n})= (k^{m-n}
\otimes_k {\mathbf 1}) \bigoplus ({\mathbf 1} \otimes_k X_{n\vert n}) $. $X_{n\vert
n}$ becomes zero in ${\cal B}_{n\vert n}$ by lemma \ref{Wakimoto}, since it is non maximal atypical. Hence
$$  \gamma(X_{m\vert n}) \ \cong\ K^{m-n}   $$
is the standard $K$-linear representation of $Gl(m-n,K)$ on
$K^{m-n}$. 
Now we use the following stronger version of the last theorem

\begin{Theorem} \label{aZ}
As a $K$-linear triangulated category the full image of each block $\Lambda$ of ${\cal T}_{m\vert n}$ in ${\cal B}_{m\vert n}$ is isomorphic to the $K$-linear triangulated category $svec_K$
spanned by the ground state $L(\lambda_0)$ of the block $\Lambda$.
\end{Theorem}

\bigskip\noindent
Theorem \ref{aZ} immediately implies theorem \ref{Z}.
Since $\gamma$ is a $K$-linear triangulated tensor functor, theorem \ref{aZ} also implies that $\gamma$ is an
exact $K$-linear equivalence of $K$-linear triangulated abelian tensor categories once we know

\begin{Lemma} \label{ss}
The category ${\cal B}_{m\vert n}$ is semisimple and hence abelian.
\end{Lemma}

\bigskip\noindent
Indeed the lemma implies exactness of the functor $\gamma$, and then corollary \ref{Wakimoto}
and theorem \ref{aZ} imply faithfulness. Hence $\varphi$ induces a faithful embedding of categories. That
$\gamma$ is full then is an immediate consequence. Put $\varphi= \gamma\circ \beta\circ \alpha$.  Then theorem \ref{aZ} and lemma \ref{ss} imply the following generalization of theorem \ref{Z}

\bigskip\noindent
{\bf Main Theorem}. {\it As a rigid $K$-linear triangulated tensor
category ${\cal B}_{m\vert n}$ is semisimple and hence abelian, and as 
a $K$-linear abelian tensor category ${\cal B}_{m\vert n}$  is equivalent to the category of $K$-algebraic finite dimensional $K$-linear
super representations of the reductive $K$-group $Gl(m-n)$. Each simple maximal atypical object $M=L(\lambda)$ maps to $$ \varphi(M) \ =\  m(\lambda) \cdot L(\lambda_0)\,[p(\lambda)]\ ,$$ where the multiplicity $m(\lambda)$ is an integer $>0$ and $L(\lambda_0)$ denotes the ground state in the block $\Lambda$ defined by $L(\lambda)$.}

\bigskip\noindent
{\it Remark}. In particular this confirms the conjecture of Kac and Wakimoto in the case of 
superlinear groups.

\bigskip\noindent
{\it Remark}. For the object $\Pi_{m\vert n} = \Lambda^{m-n}(X) \otimes Ber^{-1}$ we have
$$ \varphi(\Pi_{m\vert n}) \ =\ {\mathbf 1}_{m-n} \otimes (Ber_{n\vert n})^{-1} \ $$
in $Rep_k(H) \otimes_k {\cal B}_{n\vert n}$.
By the main theorem ${\cal B}_{m\vert n} \cong Rep_k(H) \otimes_k {\cal B}_{n\vert n}$. Hence $\Pi_{m\vert n}$ is
invertible
in the tensor category ${\cal B}_{m\vert n}$. Since $Ber_{n\vert n}$ is invertible in 
${\cal B}_{n\vert n} \cong  svec_K$ it follows that $Ber_{n\vert n} \cong {\mathbf 1}[p(Ber_{n\vert n})] = {\mathbf 1}[n]$,
hence $Ber_{n\vert n} \cong {\mathbf 1}[n]$ in 
${\cal B}_{n\vert n}$. 
Therefore theorem \ref{aZ} and lemma \ref{ss} imply

\begin{Corollary} \label{invertible}
The object $\Pi=\Pi_{m\vert n} = \Lambda^{m-n}(X) \otimes (Ber_{m\vert n})^{-1}$, where $X=k^{m\vert n}$ is the standard representation, becomes isomorphic to
${\mathbf 1}[n]$ in the triangulated tensor category ${\cal B}_{m\vert n}$
$$ \Pi_{m\vert n} \ \cong \ {\mathbf 1}[n]  \ .$$ 
\end{Corollary}

\bigskip\noindent
By lemma \ref{N-state} this in turn implies

\begin{Corollary} \label{leftshift} $L(\lambda_N)\ \cong\ L(\lambda_{N+1})[n]\ $ in $\cal B$.
\end{Corollary}

\bigskip\noindent
{\it Proof of lemma \ref{ss} using theorem \ref{aZ}}. This lemma follows from corollary \ref{abelian} of the section \ref{semisim}, since
the conditions for this corollary are provided by the parity lemma \ref{parity}
and theorem \ref{aZ}, which will be proved in the next sections \ref{moves}, \ref{three} and \ref{Proof}.

\bigskip\noindent

\bigskip\noindent
\section{\bf Basic moves} \label{moves}

\bigskip\noindent
We consider blocks $\Lambda$ for the group $Gl(m\vert n)$ of maximal atypical type. As explained in section \ref{weights} they are described by an associated set of $m-n$ crosses $\times$ on the numberline $\Z$.
The weight $\lambda$ in this block is uniquely described by $n$
labels $\vee$, which are at position different from the crosses. Attached to
a weight $\lambda$ is its cup diagram $\underline{\lambda}$ (right
move) and the oriented cup diagram $\underline{\lambda}\lambda$.

\bigskip\noindent
{\it Some simplification}. In the cup diagrams of [BS1] for many arguments the crosses $\times$ often do not play
a role. This is also true for our discussion below. Hence, for the simplicity of
exposition, we often assume $m=n$ in this section, although all statements hold for $m\geq
n$ without changes. So assume $m=n$. Then ${\cal B}={\cal B}_{n\vert n}$, so that there are no
crosses for maximal atypical weights.  The $n$ labels $\vee$ attached to a maximal atypical weight define a subset $J=\{x_1,..,x_n\}$
of the numberline $\Z$. We order the integers such that $x_1> ...
>x_n$ and  put $\lambda_j = x_j + j-1$. Then
$\lambda=(\lambda_1,..,\lambda_n;-\lambda_n,..,-\lambda_1)$ gives the
associated weight vector of a maximal atypical simple object $L(\lambda)$.

\bigskip\noindent
{\it Sectors and segments}. Every cup diagram for a weight with $n$ labels $\vee$ contains $n$
lower cups. Some of them may be nested. If we remove all inner parts of the nested cups
there remains a cup diagram defined by the (remaining)  outer cups. We enumerate these  cups from left to right. The starting points of the $j$-th lower cups is denoted $a_j$, its endpoint is denoted $b_j$. 
Then there is a label $\vee$ at the position $a_j$ and a label $\wedge$ at position $b_j$.
The interval $[a_j,b_j]$ of the numberline will be called the $j$-th sector of the cup diagram.
Adjacent sectors, i.e with $b_j=a_{j+1} -1$ will be grouped together into segments. The segments 
again define intervals in the numberline. Let $s_j$ be the starting point of the $j$-th segment
and $t_j$ the endpoint of the $j$-th segment. Between any two segments there is a distance at least $\geq 1$. The interior $I^0$ of a sector, which is obtained by removing the start and end point of the sector, always is a segment. Hence sectors, and therefore also segments have even length.

\bigskip\noindent
{\it Example $n=2$}. For the weight
$$ \mu \quad ... \ \wedge\! \wedge\! \wedge\! \wedge\! \vee\! \vee\!   \wedge\! \wedge\! \wedge\!
\wedge \ ... \ \ ,$$ with labels $\vee$ at the positions $j,j+1$ and all
other labels equal to $\wedge$, the cup diagram $\underline{\mu}$ is
described by one segment (which is a single sector)
$$ [j,j+1,j+2,j+3] \ .$$
Graphically it corresponds to a nested pair of outer cups, one
from $j+1$ to $j+2$, and one below from $j$ to $j+3$.

\bigskip\noindent
Now we fix some weight, which we denote $\lambda_{\vee\wedge}
=(\lambda_1,..,\lambda_n,-\lambda_n,..,-\lambda_1)$ for reasons to become clear
immediately. For the weight $\lambda_{\vee\wedge}$ we pick one of
the labels $x_j \in J$ at the position $i:=x_j$ such that $i+1$ is
not contained in the set of labels $J$ of the weight
$\lambda_{\vee\wedge}$. Equivalently this means $\lambda_j <
\lambda_{j+1}$ in terms of the weight vector. We define a new weight $\lambda$ (which is in
another block, and in particular is not maximal atypical) by
replacing in $\lambda_{\vee\wedge}$ the label $\vee$ at the
position $i$ by a cross $x$, and the label $\wedge$ at the
position $i+1$ by a circle $\circ$.  Attached to this new weight $\lambda$ is an
irreducible, but not maximal atypical representation $L(\lambda)$.

\bigskip\noindent
Now consider the functor $F_i$ defined in [BS4] on p.6ff and
p.10ff, which is attached to the admissible matching diagram $t$
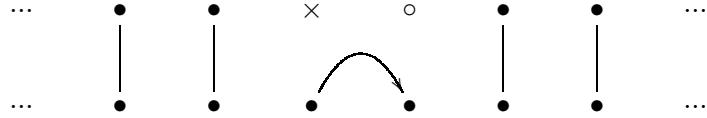
$$ \xymatrix{...& \bullet \ar@{-}[d] & \bullet\ar@{-}[d] & \times & \circ & \bullet \ar@{-}[d] &  \bullet\ar@{-}[d] & ... \cr
... & \bullet &\bullet & \bullet \ar@/^7mm/[r] & \bullet & \bullet
& \bullet &  ...\cr   }
$$ with $\times$ at position $i$ and $\circ$ at position $i+1$, and the maximal atypical object
$$ F_i(L(\lambda)) \ = \ {\bf F}_{\lambda} \ .$$
According to [BS2], lemma 4.11 this object is indecomposable and
maximal atypical with irreducible socle and cosocle
isomorphic to $L(\lambda_{\vee\wedge})$.

\bigskip\noindent
\begin{Lemma}\label{Loewy}{\it The Loewy diagram of ${\bf S}_{\lambda}$ looks
like $$ \xymatrix{ L(\lambda_{\vee\wedge}) \ar@{-}[d]\cr
{F}_{\lambda} \ar@{-}[d]\cr   L(\lambda_{\vee\wedge}) \cr}
$$ with a semisimple object ${F}_{\lambda} $ in the middle.}
\end{Lemma}

\bigskip\noindent
For the proof we give a description of the simple constituents of
${F}_{\lambda} $ below using [BS4] case (v), subcase (b), which shows
that all of these constituents have the same parity (different
from the parity of $\lambda_{\vee\wedge}$). This suffices to show the claim that $F_\lambda$ is semisimple, using lemma \ref{parity}. QED

\bigskip\noindent
Next we quote from [BS4] formula (2.13) and corollary 2.9 (of course for arbitrary $m\geq n$)

\bigskip\noindent
\begin{Lemma} {\it ${\bf F}_{\lambda}$ is a direct summand of the
representation $L(\lambda)\otimes X$, where $X$ denotes the
standard representation on $k^{m\vert n}$}. \end{Lemma}

\bigskip\noindent
Since $L(\lambda)$ is not maximal atypical, it becomes trivial in
$\cal B$ by corollary \ref{Wakimoto}. Hence the same holds for the tensor product
$L(\lambda)\otimes X$, and any of its direct summands.

\bigskip\noindent
\begin{Corollary} \label{null} {\it ${\bf F}_\lambda \cong 0$ in $\cal B$. }
\end{Corollary}

\begin{Corollary}  \label{sum} {\it $F_\lambda[1] \cong L(\lambda_{\vee\wedge}) \oplus L(\lambda_{\vee\wedge}) = 2 \cdot   L(\lambda_{\vee\wedge})$. }
\end{Corollary}

%Let ${\cal Z}_\Lambda$ be a thick triangulated subcategory of $\cal B$. 
%Then $F_\lambda \in {\cal Z}_\Lambda$ if and only if $L(\lambda_{\vee\wedge}) \in  {\cal Z}_\Lambda$.}

\bigskip\noindent
{\it Proof}. Corollary \ref{null} gives a distinguished triangle in $\cal B$
$$ L(\lambda_{\vee\wedge})[-1] \to F_\lambda \to L(\lambda_{\vee\wedge})[1] \to L(\lambda_{\vee\wedge})  \ $$
whose last arrow vanishes by lemma \ref{parity}. This proves the claim, since
$1[-1] \cong 1[1]$ holds in $\cal B$.

\bigskip\noindent
{\it The rules of [BS2], theorem 4.11}. The  constituents of ${\bf
F}_\lambda$ correspond to the maximal atypical weights $\mu$ with
defect $n$ such that \begin{enumerate} \item The (unoriented) cup
diagram  $\underline\lambda$ is a lower reduction of the oriented
cup diagram $\underline\mu t$ for our specified matching diagram
$t$. \item The rays in each "lower line" in the oriented diagram
$\underline{\mu}\mu t$ are oriented so that exactly one arrow is
$\vee$ and one arrow is $\wedge$ in each such line. \item $\mu$
appears with the multiplicity $2^{n(\mu)}$ as a constituent of
${\bf F}_\lambda$, where $n(\mu)$ is the number of "lower circles"
in $\mu t$. \end{enumerate} We remark that the lower reduction
(for more details see [BS] II, p.5ff) is obtained by removing all
"lower lines" and all "lower circles" of the diagram $\mu t$, i.e.
those which do not cross the upper horizontal numberline.

\bigskip\noindent
Let $I$ be the set of labels $\vee$ defining the maximal atypical weight
$\lambda_{\vee\wedge}$. Then $i\in I$.
To evaluate these conditions in more detail consider the 
segment $J$ of $I$ containing $i \in I$. Then also $i+1\in J$. Notice that $J$ is an interval. This
segment decomposes into a disjoint union of sectors, which
completely cover the interval $J$. We distinguish two cases

\bigskip\noindent
{\it The unencapsulated case}. Here the interval $[i,i+1]$ is one
of the sectors of $J$. We write $J=[a+1,...,i,i+1,...,b-1]$ for
the segment and call $a$ and $b$ the left and right boundary lines
of the segment. Then the label of $\lambda$ at $a$ and $b$ must be
$\wedge$ by definition. We write $I=[a,..,b]$.

\bigskip\noindent
{\it The encapsulated case}. By definition this means that the interval $[i,i+1]$ lies
nested inside one of the sectors of $J$. Hence there exists a
maximal $a < i$ defining a left starting point of a cup within the cup diagram of $\lambda$,
that has right end point $b$ such that $i+1 < b$. We write
$I=[a,...,i,i+1,...,b]$ for this subinterval of $J$ and call $a$
and $b$ the left and right boundary of $I$. The label of $\lambda$
at $a$ is $\vee$ and the label at $b$ is $\wedge$ by definition.

\bigskip\noindent
In both cases consider the sectors within $I^0=[a+1,...,b-1]$. By the maximality a $a$
$[i,i+1]$ is one of the sectors of $I^0$. 
The remaining sectors to
the left of $[i,i+1]$ and to the right of $[i,i+1]$ will be called the
lower and upper internal sectors. Let $a_j$ denote  the left starting points and
$b_j$ the right ending point of the $j$-th internal sector. The labels
at the points $a_j$ are $\vee$, and the labels at the points $b_j$ are $\wedge$.
There may be no such internal upper or lower sectors. If there are, then we will see that to each
of them corresponds an irreducible summand $L(\mu)$ of $S_\lambda$, which we will see is uniquely
described by the corresponding internal sector.

\bigskip\noindent
We summarize. In both cases the interval $I^0$ is completely
filled out by the disjoint union of the internal sectors, and one of
these internal sectors is $[i,i+1]$.

\bigskip\noindent
{\bf List of summands of ${\bf F}_\lambda$.}

\bigskip\noindent
\begin{itemize}
\item {\it Socle and cosocle}. They are defined by $L(\mu)$ for $\mu=\lambda_{\vee\wedge}$.
\item {\it The upward move}. It corresponds to the weight $\mu = \lambda_{\wedge\vee}$ which
is obtain from $\lambda_{\vee\wedge}$ by switching $\vee$ and
$\wedge$ at the places $i$ and $i+1$. It is of type
$\lambda_{\wedge\vee}$.
\item {\it The nonencapsulated boundary move}. It only occurs in the nonencapsulated
case. It moves the $\vee$ in $\lambda_{\vee\wedge}$ from position
$i$ to the left boundary position $a$. The resulting weight $\mu$ is of type
$\lambda_{\wedge\wedge}$.
\item {\it The internal upper sector moves}. For every internal upper sector
$[a_j,b_j]$ (i.e. to the right of $[i,i+1]$) there is a summand
whose weight is obtained from $\lambda_{\vee\wedge}$ by moving the
label $\vee$ at $a_j$ to the position $i+1$. These moves define new weights $\mu$ of
type $\lambda_{\vee\vee}$.
\item {\it The internal lower sector moves}. For every internal lower sector
$[a_j,b_j]$ (i.e. to the left of $[i,i+1]$) there is a summand
whose weight is obtained from $\lambda_{\vee\wedge}$ by moving the
label $\vee$ from the position $i$ to the position $b_j$. These
moves define new weights $\mu$ of type $\lambda_{\wedge\wedge}$.
\end{itemize}

\bigskip\noindent
{\it Proof of lemma \ref{Loewy}}. Except for the first in the list of summands of ${\bf F}_\lambda$, the moves of this list define  
the weights $\mu$ of the constituents $L(\mu)$ of $F_\lambda$. The parity of these weights $\mu$ is always different from $\lambda_{\vee\wedge}$. The
reason for this is, that sectors always have even length. The unique label
$\vee$ changing its position during the move, is moved by an odd number of steps.
As already explained, this suffices to prove lemma \ref{Loewy}. QED

\bigskip\noindent
{\it Remark 1}. All upper and lower internal sector moves
change the weight $\lambda_{\vee\wedge}$ into weights
$\mu$, whose cup diagram restricted to $I^0$ has a strictly smaller
number of sectors. Hence in the nonencapsulated case, the full cup
diagram of any of these $\mu$ has a strictly smaller number of sectors than
the cup diagram of $\lambda_{\vee\wedge}$.

\bigskip\noindent
{\it Remark 2}. Similarly, the nonencapsulated boundary move
changes the starting weight $\lambda_{\vee\wedge}$ into a weight
$\mu$, whose cup diagram has a strictly smaller number of sectors
except for the case $a=i-1$ (where $a=i-1$  is equivalent to $\lambda_{j-1} <
\lambda_j$).

\bigskip\noindent
{\it Remark 3}. Except for the first case in the list of summands of ${\bf F}_\lambda$, all
other moves belong to diagrams
without "lower circles". Hence $n(\mu)=1$ holds in these cases.

\bigskip\noindent
{\it Remark 4}. In the encapsulated case the diagrams $\underline{\mu} t$
do not contain "lower lines".

\bigskip\noindent
\section{Three algorithms} \label{three}

\bigskip\noindent
For $Gl(m\vert n)$
we discuss now three algorithms, which can be successively applied to a  cup diagram 
of some maximal
atypical weight within a block $\Lambda$ to
 reduce this weight to a collection of the ground state weights of this block $\Lambda$ which have the form 
$(\lambda_1,..,\lambda_{m-n},-N,...,-N;N,...,N)$ for certain large
integers $N \geq 0$.  Notice, the integers $\lambda_1,..,\lambda_{m-n}$
are fixed and describe the given block $\Lambda$. 

\bigskip\noindent
In fact, since these algorithms applies within a fixed maximal
atypical block $\Lambda$, it suffices to describe these algorithms in the case $m=n$. This 
simplifies the exposition. For this purpose assume $m=n$.

\bigskip\noindent
{\bf Algorithm I}. The first algorithm deals with a union of different
segments. The aim is to move all labels $\vee$ to the left in order
to eventually reduce everything to a single segment. For a given maximal
atypical weight $\lambda$ let $S_j=[s_j,t_j]$ from left to right
denote the segments of its cup diagram $\underline\lambda$. Let
denote $0\leq c_j =\# S_j \leq n$ their cardinalies and let denote $-\infty \leq d_j
= 1 - \vert s_{j+1} - t_i\vert \leq 0$ the negative 
distance between two neighbouring segments. We endow the set $C$ of pairs of integers
$\gamma=(c,d)$ with the lexicographic ordering. Next we endow the
set $C^n$ of all $((c_1,d_1),(c_2,d_2),.....) =
(\gamma_1,\gamma_2,.....)$ with the corresponding induced lexicographical
ordering. A cup diagram defines a maximal element in this ordering
if and only if it contains a single segment, in which case
$\gamma_1=(n,-\infty)$.

\bigskip\noindent
{\bf Claim}. Moving the starting point of the second segment to
the left increases the ordering. To be more precise: Suppose there
exist at least two segments in the cup diagram. Put $i=s_2-1$ and
$i+1=s_2$ and the weight $\lambda_{\vee\wedge}$ obtained from the
given weight $\lambda_{\wedge\vee}$ defining the cup diagram $c$
by interchange at $i$ and $i+1$. Then $[i,i+1]$ is a sector  of
the new cupdiagram $c'$ obtained in this way attached to
$\lambda_{\vee\wedge}$. Let $[a_j,b_j]$ denote the sectors of the
second segment $S_2$ with $a_1=s_1$. There are two possibilities:
\begin{itemize}
\item Then $[i,i+1][a_{1}+1,b_{1}-1]$ is a segment of $c'$ (namely
the second segment, whose first sector is $[i,i+1]$. This is the
case if and only if $d_1 < -1$;
\item or $[s_1,..,t_1][i,i+1][s_2+1,..,b_1 -1]$ is the first
segment of $c'$. This is the case if and only if $d_1=-1$.
\end{itemize}
In the first case $s'_1 =s_1$ but $d'_1 > d_1$. In the second case
$s'_1 = s_1 +1 $. Hence $c'$ is larger than $c$ with respect to
our ordering.

\bigskip\noindent
Now we consider the (unencapsulated) move centered at $[i,i+1]$ for the cup
diagram $c'_0$ attached to the weight $\lambda_{\vee\wedge}$.
Moving up gives the cup diagram $c$ we started from. The down
move, corresponding to the left boundary move, either gives as
second segment $[i-1,i]$ with unchanged first segment. Or, if
$d_1=-1$, it increases the cardinality of the first segment to
$s_1+1$. The same holds for all internal lower sector moves.
Finally for the internal upper sector moves. All these moves give
cup diagrams of the following type: With second segment
$[i,i+1][a_2+1,.., ]..[.. b_2-1]$ if $d_1 < -1$ or with first
segment $[s_1,..,t_1][i,i+1][a_2+1,.., ]..[.. b_2-1]$ if $d_1 <
-1$. Indeed they all have the same segment structure as $c'_0$,
but different sector structure. However we see that algorithm I
relates the given cup diagram $c$ to a finite number of cup
diagrams $c'$ such that $c' > c$ in our lexicographic ordering.

\bigskip\noindent
{\bf Algorithm II}. Decreasing the number of sectors within a
segment. Suppose $c$ is a maximal atypical cup diagram attaches to
a weight $\lambda_{\wedge\vee}$ with only one segment. Let
$[a_j,b_j]$ for $j=1,..,r$ denote its sectors, counted from left
to right. Assume there are at least two sectors, i.e. assume
$r>1$. Put $i=b_j$ and $i+1=a_{j+1}$ for some $1\leq j < r$. For
this recall, that any sector starts with a $\vee$ and ends with a
$\wedge$. Define $\lambda_{\vee\wedge}$ by exchanging the position
of $\vee$ and $\wedge$ in $\lambda_{\wedge\vee}$ at $i$ and $i+1$.
This gives a new cup diagram $c'_0$. It has only one segment, the
same as the segment of $c$. However the $j$-th and the $j+1$-th
sectors have become melted in $c'_0$ into one single sector. The
other sectors remain unchanged. So the numbers of sectors in the
segment decreases by one. Now consider the (encapsulated) move at
$[i,i+1]$ starting from the cup diagram $c'0$. Its move up gives
the cup diagram $c$ we started from. All internal lower and upper
moves occur within the sector $[a_j,..,b_{j+1}]$, i.e. the lower
bound $a\geq a_j$ and the upper bound is $b\leq b_{j+1}$. None of
these moves changes the cup starting from $a_j$ and ending in
$b_{j+1}$. Hence the internal moves all yield cup diagrams with
the same sector structure as $c'$. Hence algorithm II relates the
given cup diagram $c$ (with one segment $S$ and $r$ sectors) to a
finite number of cup diagrams $c'$, each of them with the same
segment $S$ but with $r-1$ sectors.

\bigskip\noindent
{\bf Algorithm III}. Now assume $c$ is a cup diagram with one
segment, which consists of a single sector $[a,...,b]$. The sector
cup from $a$ to $b$ encloses an internal cup diagram with $n-1$
labels $\vee$. This internal cup diagram necessarily defines one
segment, namely the segment $[a+1,..,b-1]$. We now apply algorithm
II to this internal segment. This finally ends up into some Kostant
weights (see [BS] II, lemma 7.2 and section \ref{Kosw})

\bigskip\noindent
{\it Further iteration}. We remark that we can start all over again and move the left starting point of the sector of a Kostant weight further to the left
using algorithm I, and then repeat the whole procedure of applying algorithms I, II and III. At the end this allows to replace the given Kostant weight by some other Kostant weights further shifted to the left on the numberline (with all crosses $\times$ removed in case $m\geq n$). If we repeat this down shift of Kostant weights sufficiently often we end up with a bunch of Kostant weights, that are 
ground states of the block, i.e. whose associated irreducible representation is one of the ground state
representations $L(\lambda_N)$ for large $N$.

\bigskip\noindent

\section{Proof of theorem \ref{aZ}} \label{Proof}

\bigskip\noindent
To prove the theorem \ref{aZ} we now fix a maximal atypical block $\Lambda$ of $\cal T$
and its ground state representation $L=L(\lambda_0)$.
Consider the thick triangulated subcategory ${\cal Z}_{\Lambda}$ of $\cal B$ associated to $L$ as defined in section \ref{groundcat}.
 To show that a given simple maximal atypical representation $L(\mu)$ of $\Lambda$ has image in ${\cal Z}_\Lambda $
 it suffices that it is zero in ${\cal B}/{\cal Z}_\Lambda$. If this holds for all simple objects of the block $\Lambda$, then it also holds for all objects of the block $\Lambda$. 
 
\bigskip\noindent
An object $A$ will be called a virtual ground state object, if there exists an isomorphism
in $\cal B$ of the form $A \oplus A' \cong A''$ where $A'$ and $A''$ is isomorphic
to a finite direct sum of higher ground state objects $L_N$ of the block $\Lambda$.
We can apply the algorithms I, II and III and corollary \ref{sum} to show by induction that
there exist virtual ground state objects $Y$ and $Y'$ and an isomorphism in $\cal B$ 
$$  L(\mu) \oplus Y \ \cong \ Y' \ .$$
This immediately implies
that also $L(\mu)$ is a virtual ground state object. 

\bigskip\noindent
In lemma \ref{down} we show,
that all higher ground states $L_N$ of $\Lambda$ (for all $N\geq 0$)
are in ${\cal Z}_\Lambda$. 
For this we use the next 

\bigskip\noindent
{\bf Algorithm IV}. Let $\lambda$ be a Kostant weight. By [BS] II,
lemma 7.2 this means that the labels $\vee$ of $\lambda$ define an
interval $[a,..,a+n-1]$ after the crosses have been removed. Starting from the label
$\lambda$ we make successively moves with the right most label $\vee$ away
to the right by $i$ steps. Let us call these new weights $S^i$ so that
$S^0$ is the Kostant weight we started from, and so on. Put
$i=a+n-1$. Make a first move with $[i,i+1]=[a+n-1,a+n]$. This move
is encapsulated and gives the Loewy diagram
%$$ \xymatrix{ L(\lambda_{\vee\wedge}) \cr
%L(\lambda_{\wedge\vee})\ar@{-}[u]\cr
%L(\lambda_{\vee\wedge})\ar@{-}[u] \cr}
%$$
%or
$$ \xymatrix@C-9mm{ S^0 \cr
S^1\ar@{-}[u]\cr S^0 \ar@{-}[u] \cr}
$$
Next move for $[i+1,i+2]$ gives the Loewy diagram with four irreducible constituents
$$ \xymatrix{ & S^1 & \cr
S^2 & \oplus \ar@{-}[u]  & S^0 \cr & S^1\ar@{-}[u] & \cr}
$$
and so one until the first move, which is not encapsulated. Here
we end up in a Loewy diagram of type
$$ \xymatrix@C-9mm{ & S^{n-1} & \cr
S^n \ \ \ \oplus & \Pi \ar@{-}[u]  &  \oplus \ \ \ S^{n-2}\cr &
S^{n-1}\ar@{-}[u] & \cr}
$$
with additional fifth constituent $\Pi$, where $\Pi$ is of Kostant type in the given block such that compared to the Kostant weight $\lambda$ we started from all labels
$\vee$ have been shifted to the left by one, and hence are at the positions $[a-1,...,a+n-2]$.

\bigskip\noindent
\begin{Lemma} \label{down}{\it Suppose $L(\lambda_N) \in {\cal Z}_\Lambda$, then 
also $L(\lambda_{N+1})\in {\cal Z}_\Lambda$.}
\end{Lemma}

\bigskip\noindent
{\it Proof}. Use that ${\cal Z}_\Lambda$ is a thick triangulated subcategory of $\cal B$
by lemma \ref{thick}. Hence it suffices that $S^0 =0$ in ${\cal B}/{\cal Z}_\Lambda$ implies
$S^i=0$ and hence $\Pi=0$ in ${\cal B}/{\cal Z}_\Lambda$, which obviously follows from
the Loewy diagrams displayed above. Since $\Pi = L(\lambda_{N+1})$, if $S^0=L(\lambda_N)$, we are done. QED

\bigskip\noindent
This shows that for every simple object $L(\mu)$ there exist $A$ and $A'$ in
${\cal Z}_\Lambda$ such that $L(\mu) \oplus A \cong A'$ in $\cal B$.
Hence $L(\mu) =0$ in ${\cal B}/{\cal Z}_\Lambda$. Therefore $L(\mu)\in {\cal Z}_\Lambda$.
By parity reasons therefore $L(\mu) \cong m(\mu) \cdot L[p(\mu)]$ for the uniquely defined
parity $p(\mu)$ (which is easily computed from the Bruhat distance from the ground state $\lambda_0$). This proves theorem \ref{aZ}.

\bigskip\noindent
%{\it Proof}. Suppose we would like to known wether $L(\mu)\in \cal
%Z$. We can assume that $\mu$ is not basic. Pick some $i$ lying
%between two segment such that $i+1$ is in the segment.
% Suppose
%$\lambda_{\vee\wedge}$ is obtained from $\mu=\lambda_{\wedge\vee}$
%by an exchange at the positions $i$ and $i+1$. Then we claim
%$$ (*) \quad \quad \lambda_{\vee\wedge}\in \cal Z \Longrightarrow \lambda_{\wedge\vee} \in
%{\cal Z} \ .$$ Indeed, since ${\bf F}_\lambda$ is zero in $\cal
%B$, we conclude from the last lemma (semisimplicity of
%$F_\lambda$) and $\lambda_{\vee\wedge}\in \cal Z$ that all
%summands of $F_\lambda$ are in $\cal Z$. This proves our claim.

%\bigskip\noindent
%We use this to prove the lemma by induction. Let $c_1$ be the
%cardinality of the bottom segment of a weight, $c_2$ the
%cardinality of the second segment from left, and so on. This
%attaches a sequence of numbers $c_1,c_2,..$ to each weight $\mu$.
%Among the weights $\mu$ for which $L(\mu)$ is not in $\cal Z$
%choose $\mu$ such that this associated sequence of numbers
%$c_1,...$ is maximal in the lexicographic ordering. Since the
%passage from $\mu=\lambda_{\wedge\vee}$ to
%$\lambda=\lambda_{\vee\wedge}$ never increases this lexicographic
%ordering, we can modify $\mu$ by moving certain $\vee$ in the cup
%diagram of $\mu$ to the left. Although the lexicographic ordering
%may not increase for a single move of this type, it does so after
%finitely many steps as long as $\mu$ is not basic. Hence, using
%the implication (*) sufficiently often, we obtain a contradiction.
%QED

\section{Odds and ends}

\bigskip\noindent
Consider a maximal atypical block of $\cal T$.
Since $Ber_{m\vert n}$ is invertible, the tensor product with
$Ber_{m\vert n}$ defines equivalences between maximal atypical blocks
and their twisted images. Hence using a twist by a power of the Berezin we may, without restriction of generality, assume 
that the block $\Lambda$ contains a ground state weight vector of the special form $\lambda_0 =(\lambda_1,...,\lambda_{m-n}\ ,\ 0,...,0\ ;\ 0,...,0)$ where
$\lambda_1 \geq \lambda_2 \geq ... \geq \lambda_{m-n}=0$. For this see the remarks following lemma \ref{ground}. 
Hence we may assume that the ground state is a covariant representation $L\cong \{ \lambda \}$
associated to the partition $\lambda_1 + \lambda_2 + .. + \lambda_{m-n-1}$. 
We say that $\Lambda$ is a block with normalized ground state.

\bigskip\noindent
Consider the $K$-linear triangulated tensor functor
$$ \varphi: {\cal T}_{m\vert n} \longrightarrow Rep_k(H) \otimes_k {\cal B}_{n\vert n} \ $$
for $L=\{ \lambda \} = Schur_\lambda(X)$. Since $\varphi(X) = k^{m-n} \otimes_k {\mathbf 1}$ (the standard representation is not maximal atypical for $m=n$ and corollary \ref{Wakimoto}),
we conclude $$\varphi(L) \ \cong\ Schur_\lambda(k^{m-n}) \otimes_k {\mathbf 1}\ .$$
Since $\mathbf 1$ is the ground state of the unique maximal atypical block of ${\cal T}_{n\vert n}$,
this implies

\begin{Lemma} \label{groundtoground} For blocks with normalized ground states
the functor $\varphi$ maps ground states to ground states. 
\end{Lemma}

\begin{Lemma} \label{Ber} $\varphi(Ber_{m\vert n}) = det \otimes_k Ber_{n\vert n} =det \otimes {\mathbf 1}[n]$. 
\end{Lemma}

\bigskip\noindent
{\it Proof}. Obvious.

\bigskip\noindent

\section{Multiplicities}\label{Weyl}

\bigskip\noindent
Fix a maximal atypical block $\Lambda$. Let the ground state vector
of $\Lambda$ be $$\lambda = (\lambda_1,...,\lambda_{n-m},M,..,M;-M,...,-M)$$ for
$M=\lambda_{n-m}$. The block $\Lambda$ is characterized by $(\lambda_1,..,\lambda_{m-n})$
respectively the corresponding irreducible representation $\rho= Schur_{\lambda_1,..,\lambda_{m-n}}(k^{m-n})$ in $Rep_k(H)$ with the following convention. For $M=\lambda_{m-n} < 0$ we
define $Schur_{\lambda_1,..,\lambda_{m-n}}(k^{m-n}) := Schur_{\lambda_1-M,..,\lambda_{m-n}-M}(k^{m-n}) \otimes det^M $ by abuse of notation. By Lemma \ref{Ber} this notation
behaves nicely with respect to the triangulated tensor
functor 
$$ \varphi: {\cal T}_{m\vert n}\ \longrightarrow\ Rep_k(H) \otimes_k {\cal B}_{n\vert n} \ .$$
Indeed, $\varphi(L(\lambda)) = \varphi( Ber_{m\vert n}^M \otimes Schur_{\lambda_1-M,..,\lambda_{m-n}-M}(X)) $ for $X=k^{m\vert n}$ coincides with 
$ det^{M} \otimes Schur_{\lambda_1-M,..,\lambda_{m-n}-M})(\varphi(X)) \otimes_k Ber_{n\vert n}^M$. Since $\varphi(X) \cong k^{m-n} \otimes_k 1$ in $\cal B$, this gives with the convention above  $\varphi(L(\lambda)) =  Schur_{\lambda_1,..,\lambda_{m-n}}(k^{m-n}) \otimes_k Ber_{n\vert n}^M$.
Recall ${\cal B}_{n\vert n} = {\cal Z}\cong svec_K$. Using corollary \ref{leftshift}
an obvious reexamination of the proof of theorem \ref{aZ} 
now shows 

\bigskip\noindent
\begin{Theorem} \label{indep} {\it For each weight $\mu$ in the fixed maximal atypical block $\Lambda$
we have 
$$  \varphi(L(\mu)) \ \cong \ m(\mu) \cdot Schur_{\lambda_1,..,\lambda_{m-n}}(k^{m-n}) \otimes_k  {\mathbf 1}[p(\mu)] $$
in $Rep_k(H) \otimes {\cal Z}$ for some integral multiplicity $m(\mu) \geq 1$, which
only depends on the relative position of the weight $\mu$ with respect to
the ground state weight, considered on the numberline $\Z$
with all crosses $\times$ removed. For the ground states the multiplicity is one.}\end{Theorem}

\bigskip\noindent
In particular, this theorem shows that the computation of the multiplicities $m(\mu)$ can be reduced to the case $m\! =\! n$. The computation of the parity $p(\mu)=\sum_{i=1}^n \mu_{m+i}$
is reduced to the case of the ground state. By lemma \ref{Ber} one can reduce to the case of a block with a normalized ground state, where the parity is even by the computation preceding
the theorem. 

\bigskip\noindent
{\it The muliplicities $m(\mu)$}. As already explained, as a consequence of theorem \ref{indep},
we may assume $m\! =\! n$ for the computation of the multiplicities. These multiplicities 
are numbers attached to cup diagrams $\underline\mu$ with $n$ cups (and without lines). We have already shown that the multiplicity $m(\mu)$ is one for the ground state $\mu=\mathbf 1$.
The same holds for all powers $Ber^N$ of the Berezin by lemma \ref{Ber}. Hence
for any completely nested cup the multiplicity $m(\mu)$ is one. To deal with a general maximal atypical weight
our strategy is the following. We consider cup diagrams for various $n$ with the aim to reduce the computation of $m(\mu)$ for a cup diagram with $n$ cups to the case of cup diagrams with $<n$
cusps. 

\bigskip\noindent
For a completely nested cup the multiplicity $m(\mu)$ is one.
In general let $\mu$ have the sectors $S_1,..,S_r$ with length $2n_1,..,2n_r$ with corresponding partial cup diagrams
$\underline{\mu_1},...,\underline{\mu_r}$. Notice $n=n_1+...+n_r$. Each $S_i$ defines a number interval $[a_i,b_i]$. Using the Berezin we see that the multiplicity does not change under
a translation of the cup diagram.
Now the algorithms II and III applied to the cup diagram of $\mu$ show, that all nested cups $\underline{\mu_i}$ can be reordered to become completely nested without destroying the sector structure of the original cup diagram $\mu$.
In this way the cup diagram can be rearranged so that all nested cups are completely nested
cups. This process proves the formula
$$ (*) \quad \quad  m(\mu) \ = \ m(\nu) \cdot \prod_{i=1}^r \ m(\mu_i) \ ,$$
where $\nu$ is the cup diagram with the same sectors as $\mu$, but so that each sector $S_i$
defines a maximal nested cup diagram with labels $\vee$ at the position $a_i,a_i+1,...,a_i + n_i$.

\bigskip\noindent
\begin{Lemma} \label{mult} Suppose $\nu$ is a maximal atypical weight $\nu$ with $r$ sectors. If all sectors of the cup diagram of $\nu$ have completely nested cup diagrams of lengths say $2n_1,..,2n_r$, then
$$m(\nu) = { n \choose n_1,\ldots ,n_r } \quad  \mbox{(multinomial coefficient)} \ .$$
\end{Lemma}

\begin{Lemma}  For maximal atypical weights $\mu$ with $n$ labels $\vee$
the multiplicity  
$$  m(\mu) \ = \  { n \choose n_1,\ldots,n_r } \cdot\ \prod_{i=1}^r \ m(\mu_i) $$
satisfies the inequality $$1 \leq m(\mu) \leq n! \ .$$ Equality
at the right side holds if and only if the cup diagram of $\mu$ is completely unnested (i.e. all sectors have length 2).  Equality on the left holds if and only if the cup diagram has only one sector which is a completely nested sector (translates of the ground state). \end{Lemma}

\bigskip\noindent
{\it Proofs}. Induction on $n$ using formula (*), theorem \ref{indep}
and lemma \ref{mult}.

\bigskip\noindent

\section{Proof of lemma \ref{mult}}

\bigskip\noindent
{\it Case of one segment}. Suppose $\nu$ has only one segment.
Then $m(\nu)=m(n_1,..,n_r)$ depends only on the sector lengths $n_1,..,n_r$.
Algorithm II applied to the first two sectors, combined with formula (*) from above, gives the recursion formula
$$ 2 \cdot m(u\! +\! v, n_3,..,n_r) m(u\! -\! 1,1,v\! -\! 1) \ - \ m(u,v,n_3,..,n_r) \ = \ $$
$$ m(u\! +\! v,n_3,..,n_r)m(u\! -\! 1,v)m(1,v\! -\! 2)\
+\ m(u\! +\! v,n_3,..,n_r) m(u,v\! -\! 1) m(u\! -\! 2,1) $$  
for $m(u,v,n_3,..,n_r)$ in $u$ and $v$. All terms except $m(u,v,n_3,..,n_r)$ involve 
either fewer variables or less labels $\vee$.  
This allows to verify the claim by induction on the number $\sum_{i=1}^r n_i$ of labels $\vee$ and then of sectors $r$. The verification of the induction start $r=1$ is obvious by definition. 
 So it suffices that the multinomial coefficient $m(n_1,..,n_r)=(\sum_{i=1}^r n_i)!/\prod_{i=1}^r (n_i)!$ satisfies the recursion relation of algorithm I. The trivial property  $$m(n_1,n_2,n_3,..,n_r)\ =\ m(n_1,n_2) \cdot m(n_1+n_2,n_3,..,n_r)$$ of multinomial coefficients allows to assume $r=2$.  
The recursion formula then boils down to the identity $2uv= (u+v) + v(u-1) + u(v-1)$. This proves the assertion if there is only one segment.

\bigskip\noindent
{\it The case of more than one segment}.
Now suppose $\nu$ is a maximal atypical totally nested weight 
with $s >1$ 
segments and with a total number of $r$ sectors of lengths $2n_1,..,2n_r$. 
Notice that
all segments are sectors by our assumption on $\nu$. We then symbolically write
$$ \nu \ = \  ...\ \wedge\!\wedge\ (S_1 \!\wedge ... \ q\ ... \wedge\! S_2)\ \wedge ...\ \  \mbox{rest with higher segments}  \   $$
for the segment diagram, where $q$ denotes the distance between the first and second segment. 
To show that the multiplicity formula of lemma \ref{mult} also holds in general
we now use algorithm I to increase the size of the first sector. We assume by induction that
the formula holds for maximal atypical totally nested weight with $<s$ segments
or for maximal atypical totally nested weight with $\geq s$ segments and more than $2n_1$ elements in the first sector or  with $\geq s$ segments and $2n_1$ elements in the first sector 
but smaller distance $q$ between the first and second sector.
This start of the induction is the case with one segment already considered.

\bigskip\noindent
{\it First case}. Suppose the distance $q=1$. 

\bigskip\noindent
a) Then for completely nested sectors $S_1$ and $S_2$ of length $2n_1$ and $2n_2$
$$ \nu \ = \ ... \  \wedge (\wedge \ S_1 \wedge S_2) \wedge ...\ \  \mbox{rest with higher segments}  \ .  $$
b)  Let $\lambda_{\vee\wedge}$ denote the weight obtained by moving the starting
point of the second sector $S_2$ one step down so that that it touches the end of the first sector $S_1$. This new weight $\lambda_{\vee\wedge}$ has $s-1$ segments with segment structure 
$$   \lambda_{\vee\wedge} \ = \ ... \ \wedge (\wedge\ T_1\ \wedge) \wedge ...\ \ \mbox{rest with higher segments} $$
whose first segment $T_1$ has length $2(n_1+n_2)$ with three completely nested sectors of lengths $2n_1,2,2(n_2-1)$ respectively. 
Algorithm I gives three further weights:

\bigskip\noindent
c) The boundary move weight  with segment diagram
$$ ... \ \wedge (S'_1 \wedge S'_2 \ \wedge) \wedge ...\ \ \mbox{rest with higher segments}  \ .  $$
where the first and second segments $S'_1$ and $S'_2$ are completely nested sectors of lengths $2(n_1+1)$ and $2(n_2-1)$.

\bigskip\noindent
d) The interval lower sector move gives a weight with $s-1$ segments and diagram
$$   \ ...\  \wedge (\wedge\ T'_1\ \wedge) \wedge ...\ \ \mbox{rest with higher segments} $$
where the segment $T'_1$ of length $2(n_1+n_2)$ has two sectors. The second sector is completely nested
of length $2(n_2-1)$. The first sector $I$ has length $2(n_1+1)$ and its interior segment $I^0$ decomposes into
two completely nested sectors of lengths $2(n_1-1)$ and $2$ respectively.

\bigskip\noindent
e) The interval upper sector move gives a weight with $s-1$ segments and diagram
$$ ...\ \wedge (\wedge\ T''_1\ \wedge) \wedge ...\ \ \mbox{rest with higher segments} $$
where the first segment $T''_1$ has length $2(n_1+n_2)$ with two sectors. The first sector is completely nested
of length $2n_1$. The second sector $I$ has length $2n_2$ and its interior segment $I^0$ decomposes into
two completely nested sectors of lengths $2$ and $2(n_1-2)$ respectively.

\bigskip\noindent
Again we show that the multinomial coefficient $m(n_1,..,n_r)$ satisfies the
recursion relation of algorithm I. This suffices to prove our assertions. Again the trivial property $m(n_1,n_2,n_3,..,n_r)\ =\ m(n_1,n_2) \cdot
m(n_1+n_2,n_3,..,n_r)$ of multinomial coefficients allows to assume $r\! =\! 2$.  
The desired recursion equation of algorithm I, that the sum of the multiplicities of a), c), d) and e) is twice the
multiplicity of b), then boils down to the binomial identity
$$2 \!\cdot\! { n_1 + n_2 \choose n_1,1,n_2\! -\! 1 } = { n_1\! +\! n_2 \choose n_1 } + { n_1\! +\! n_2 \choose n_1\! +\! 1 } + { n_1\! +\! n_2 \choose n_1\! +\!1 } { n_1  \choose 1 }
+ { n_1\! +\! n_2 \choose n_1 } { n_2\! -\! 1 \choose 1}  \ .$$

\bigskip\noindent
{\it Second case}.  Now suppose $q\geq 2$ for the distance $q$ between the first and the second sector.

\bigskip\noindent
a) Then $ \nu \ = \ ...\   \wedge (S_1\! \wedge ... \ q\ ... \wedge\! S_2) \wedge ...\ \  \mbox{rest with higher segments}  $. 
The first and second segments $S_1,S_2$ are completely nested sectors of length $2n_1$ respectively $2n_2$.

\bigskip\noindent
b) Consider the weight $\lambda_{\vee\wedge}$ which is obtained by moving the starting
point of the second sector $S_2$ one step down to the left. Since $q>1$ it does not touch the end of the first sector $S_1$. The new weight $\lambda_{\vee\wedge}$ still has $s$ segments, but now with the segment structure 
$$   \lambda_{\vee\wedge} \ = \ ...\ \wedge (S_1\! \wedge ...\ q-1\ ... \wedge\! S'_2\ \wedge ) \wedge ... \ \ \mbox{rest with higher segments} $$
where the second segment $S'_2$ of lenght $2n_2$ has two completely nested sectors
of lengths $2$ and $2(n_2-1)$.
The algorithm I gives two further weights:

\bigskip\noindent
c) If $q>2$, the boundary move weight gives a diagram with $s+1$ segments 
$$ ...\ \wedge (S_1 \!\wedge ...\ q-2\ ... \wedge\! T_{23}\ \wedge ) \wedge ...\ \ \mbox{rest with higher segments}  \  $$
where $T_{23} = [\vee\wedge]\wedge T'_2$ has two completely nested segments $T'_1=[\vee\wedge]$ and $T'_2$ of lengths
$2$ respectively $2(n_2-1)$. 

\bigskip\noindent
If $q=2$ we get a diagram with $s$ segments
$$ ...\ \wedge (S'_1 \wedge T'_2 \ \wedge) \wedge ...\ \ \mbox{rest with higher segments}  \   $$
where the first segment $S'_1 $ has two completely nested sectors of lengths
$2n_1$ respectively $2$ and the second segment $T'_2$ is a completely nested sector of length $2(n_2-1)$; 
with distance $q=1$ from the first sector. 

\bigskip\noindent
e) The interval upper sector move gives a weight with $s$ segments and diagram
$$ ...\ \wedge (S_1 \!\wedge ...\ q-1\ ...\wedge\! T_2\ \wedge) \wedge ...\ \ \mbox{rest with higher segments} $$
where the second segment $T_2$ is a sector of length $2n_2$, its interior decomposes into two completely nested sectors
of lengths $2$ and $2(n_2-2)$. 

\bigskip\noindent
The recursion relation of algorithm II, that twice the 
multiplicity of b) is the sum of the multiplicities of a), c) and e), holds for
the multinomial coefficient. This amounts to
the binomial identity
$$ 2\!\cdot\! {n_1\! +\! n_2 \choose  n_1,1,n_2\! -\! 1} = 
{n_1\! +\! n_2 \choose  n_1} + {n_1\! +\! n_2 \choose  n_1} {n_2\! -\! 1 \choose 1}
+ {n_1\! +\! n_2 \choose  n_1,1,n_2\! -\! 1} \ .$$
This finally completes the proof of  lemma \ref{mult} using induction on the distance $q$.

\bigskip\noindent
\section{Appendix: The class $\xi_n$}
\label{class}

\bigskip\noindent
 Consider the exact BGG complex [BS2], thm. 7.3 for the Kostant weight $\mu=0$ given by
$ ... \to V_{2} \to V_{1} \to V_0 \to {\mathbf 1} \to 0 $
with
$$ V_{j} \ = \ \bigoplus_{\lambda \leq 1, \ l(\lambda,1)=j} \ V(\lambda) \ .$$
%Hence there is a
%spectral sequence
%$ Ext_{\cal T}^i(V_{j}, X)\ \Rightarrow\ Ext_{\cal T}^{i-j}(1, X)$.

%\bigskip\noindent
%{\it Fact 2}.  By [BS3] for all weights $\nu$
%$$ \sum_{i\geq 0} dim_k( Ext^i_{\cal T}(V(\nu),1) \ \leq \ 1 \ .$$
%{\it Fact 3}. $Ext_{\cal T}^i(V(\lambda), L(\mu))=0$ unless $\lambda\leq \mu$ and
%$i\leq l(\lambda,\mu)$. See [BS3] formula (5.3) and corollary 5.5. Furthermore
%$Ext_{\cal T}^i(V(\lambda), L(\mu))\cong k$ for $i=l(\lambda,\mu)$ and $\lambda\leq \mu$.

%

%

%\bigskip\noindent
%Suppose there exists a nontrivial morphism $\xi: 1 \to {\cal L}[n]$ in ${\cal K}$ for some power $\cal L\in \cal T$ of the Berezin $Ber$, which induces an isomorphism in $\cal B$.
%If $\cal L$ is nontrivial, then this determines $\cal L$ uniquely up to inversion (use restriction
%to $psl(n\vert n)$ and the existence of two additional morphisms $\xi_n$ and $\eta_n$ [BKN1]. Indeed, there is a second morphism $\eta: 1 \to {\cal L}^{-1}[n]$ in $\cal T$ defined
%by $\xi^*: {\cal L}[-n] \to 1 $.

\bigskip\noindent
\begin{Proposition}\label{xi} {\it For $Gl(n\vert n)$ there is a nontrivial morphism $\xi_n: {\mathbf 1}\to Ber_{n\vert n}[n]$
in ${\cal K}$ which becomes an isomorphism in $\cal B$.
Hence $Ber$ is contained in $\cal Z$.} \end{Proposition}

\bigskip\noindent
{\it Proof}. Appling the antiinvolution ${}^*$ we get
$$ 0\to {\mathbf 1} \to V_0^*\to V_1^* \to \cdots \to V_{n-1}^* \to V_n^* \to \   ,$$ which defines
a Yoneda extension class $\xi \in Ext^n(Q,1)$
$$ 0\to {\mathbf 1} \to V_0^*\to V_1^* \to \cdots \to V_{n-1}^* \to Q \to 0   $$
for $ Q \ \cong\ im(d^*: V_{n-1}^* \to V_n^*) \  \hookrightarrow \  V_n^*$.

\bigskip\noindent
Now $V_i^*=0$ in $\cal H$, since $V([\lambda])^*=0$ holds in ${\cal H}$ for all cell modules $V([\lambda]$ by the definition of $\cal H$. Hence in $\cal H$, and therefore in $\cal B$, we get
 $$Q \cong {\mathbf 1}[n] \ .$$ 
We will now construct a map $i: Ber^{-1} \hookrightarrow Q$ in $\cal T$,
which defines a nontrivial morphism in $\cal B$. Then $\xi$ is nontrivial in ${\cal K}$,
hence $i^*(\xi)$ defines a nontrivial extension
in $Ext^n_{\cal T}(Ber^{-1},1)$.

\bigskip\noindent
{\it Nontriviality}. To show that a given morphism $i: Ber^{-1} \hookrightarrow Q$ is nontrivial in $\cal B$
is equivalent to show that the transposed morphism $i^*: Q^* \to Ber^{-1}$ is nontrivial in $\cal B$. We first show that $i^*: Q^* \to Ber^{-1}$ is nonzero in $\cal H$. Then $i^*$ remains nonzero in $\cal B$, using $Q^* \cong {\mathbf 1}[-n]\cong {\mathbf 1}[n]$ in $\cal H$ and using that after restriction to $psl(n,n)$ the morphism
$i^*$ is in the central graded ring $R^\bullet_{\cal K}$ and of positive degree, hence by proposition \ref{grad} the morphism $i^*$ can not become a zero divisor for the localization $\cal B$ and thus is a nonzero morphism in $\cal B$.

\bigskip\noindent
To show that $i^*\neq 0$ in $\cal H$ we argue as follows: If $i^*=0$ in $\cal H$, then the composite morphism $V_n \twoheadrightarrow Q^* \twoheadrightarrow Ber^{-1}$, defined in $\cal T$, becomes zero in $\cal H$. Since $Ber^{-1}$ is simple and $V_n$ is a cell object, we can apply theorem \ref{hot} to obtain $Hom_{\cal H}(V_n,Ber^{-1}) = Hom_{\cal K}(V_n,Ber^{-1})=
Hom_{\cal T}(V_n,Ber^{-1})$.  Since $V_n \twoheadrightarrow Ber^{-1}$ is an epimorphism in $\cal T$, the composite map is nonzero in $\cal T$ and therefore nonzero in $\cal H$. This completes the proof that $i^*$ is not zero in $\cal H$. Therefore, once we have constructed an epimorphism  $i^*: Q^* \to Ber^{-1}$
in $\cal T$, this proves the proposition.

\bigskip\noindent
{\it Existence}.  To define $i^*$ recall, that the boundary morphisms $d^*$ are dual to the morphism $d$ defined in  [BS2]: With respect to the decomposition
$$ V_n \ = \bigoplus_{\mu \leq 0\ ,\ l(\mu,0)=n} V(\mu) \ ,$$
the $d:V_{n-1} \to V_n$ are defined as the sum of morphisms $f_{\lambda\mu}$. See [BS2], lemma 7.1. These $f_{\lambda\mu}$ are obtained as follows:
There exists a projective $P=P(\lambda)$ and an endomorphism $f:P\to P$, an filtration
$N \subset M \subset P$ with $P/M =V(\lambda)$ and $M/N = V(\mu)$ such that
$f(P)\subset M$ and $f(M)\subset N$ so that $f$ induces the morphism $f_{\lambda\mu}:P/M=V(\lambda) \to M/N=V(\mu)$. 

\bigskip\noindent
To define an epimorphism
$$ i^*: Q^* \cong Im\bigl(d:V_n \to V_{n-1}\bigr) \ \twoheadrightarrow Ber^{-1} \ $$
notice that $d:V_n \to V_{n-1}$ is $\sum d_\mu$ for $d_\mu=d\vert_{V(\mu)}: \ V(\mu) \to V_{n-1}$
where $\mu < 0$ and $l(\mu,0)=n$. The cosocle of $V(\mu)$ is the simple
object $L(\mu)$. Hence $V(\mu)=V(Ber^{-1})$ is the unique summand of $V_n$
with cosocle $Ber^{-1}$. Since none of the morphisms $d_\mu$ is trivial in $\cal T$ (see [BS2]),
the summand $Ber^{-1}$ in the cosocle of $V_n$ maps nontrivially to the cosocle
of its image $Q^*$ in $V_{n-1}$.
The only indecomposable summand $V(\mu)$ of $V_n$ containing $Ber^{-1}$ in its cosocle is $V(Ber^{-1})$. The cosocle of $V(Ber^{-1})$ therefore injects into $Q$ by the definition of the morphism $d$. This completes the proof of proposition \ref{xi}.

\bigskip\noindent

\section{Appendix: Semisimplicity}\label{semisim}

\bigskip\noindent

\bigskip\noindent
We consider a triangulated category\footnote{As pointed out by Heidersdorf, that this is a special case of a triangulated category $\cal B$ with a cluster tilting cotorsion pair $(\cal U,V)$ (see [Na]).}
$\cal B$, such that
there are strictly full additive subcategories ${\cal B}_0$ and ${\cal B}_1$ with the following properties
\begin{enumerate}
\item ${\cal B} = {\cal B}_0 \oplus {\cal B}_1 $
\item ${\cal B}_0[1] = {\cal B}_1$
\item ${\cal B}_1[1] = {\cal B}_0$
\item $Hom_{\cal B}({\cal B}_0,{\cal B}_1)=0$.
\end{enumerate}

\bigskip\noindent
By property 1. and 4. it easily follows, that the decomposition
1. is functorial such that there are adjoint functors $\tau_0:{\cal B} \to {\cal B}_0$
and $\tau_1: {\cal B\to \cal B}_1$ with functorial distinguished triangles
$(\tau_0(X),X,\tau_1(X))$, s.t. $Hom_{\cal B}(A,X)= Hom_{\cal B}(A,\tau_0(X))$
for $A\in {\cal B}_0$ and similarly $Hom_{\cal B}(A,X)= Hom_{\cal B}(A,\tau_1(X))$
for $A\in {\cal B}_1$. The next lemma immediately follows from the long exact $Hom$-sequences
attached to distinguished triangles

\bigskip\noindent
{\bf Extension Lemma}. {\it If $A,C\in {\cal B}_0$ and $(A,B,C)$ is a distinguished
triangle, then $B\in {\cal B}_0$.}

\bigskip\noindent
\begin{Lemma} {\it ${\cal B}_0$ is an abelian category.}
\end{Lemma}

\bigskip\noindent
{\it Proof}. For a morphism $f:X\to Y$ in ${\cal B}_0$ let $Z_f$ be a cone in $\cal B$.
Put $Ker_f = \tau_0(Z_f[-1])$ and $Ker_f = \tau_0(Z_f)$. Then $Ker_f,Koker_f$ are in ${\cal B}_0$
and represent the kernel resp. kokernel of $f$ in ${\cal B}_0$. This is an immediate
consequence of the long exact $Hom$-sequences and property 4. Furthermore
$=b\circ a$ factorizes in the form $a:X\to Z$ and $b:Z\to Y$ by the
octaeder axiom, where $Z$ is a cone of the composed morphism $i$ defined by
$Ker_f \to Z_f[-1] \to X$. By the octaeder axiom there exists a distinguished
triangle $(Koker_f,Z,Y)$. Hence by the extension lemma $Z\in {\cal B}_0$.
The octaeder axiom provides the morphisms $a$ and $b$ and proves
$a_*: Z \cong Koker_i$  and similarly $b_*: Z \cong Ker_\pi$ for the 
morphism $\pi: Y \to Koker_f$. Since obviously ${\cal B}_0$ is an additive subcategory
by the functoriality of the decomposition 1., this implies that ${\cal B}_0$ is an abelian 
category. QED

\bigskip\noindent
By construction the exact sequences $0\to A\to B\to C\to 0$
in ${\cal B}_0$ correspond to the
distinguished triangles $(A,B,C)$ in ${\cal B}_0$.

\begin{Lemma} \label{sem}{\it The abelian category ${\cal B}_0$ is semisimple.} 
\end{Lemma}

\bigskip\noindent
{\it Proof}. For a short exact sequence the corresponding
triangle $(A,B,C)$ splits, since the morphism $C\to A[1]$ 
vanishes by property 4. QED

\begin{Corollary}\label{abelian} {\it ${\cal B}$ is a semisimple abelian category.} 
\end{Corollary}

\bigskip\noindent
Suppose $(X,Y,Z)$ is a distinguished triangle in $\cal B$.
Then for $Z\in {\cal B}_0$ there exists an exact sequence
$$ 0 \longrightarrow \tau_0(X) \longrightarrow \tau_0(Y) \longrightarrow \tau_0(Z) $$
in ${\cal H}_0$. Similary if $X\in {\cal B}_0$ there exists
an exacts sequence
$$  \tau_0(X) \longrightarrow \tau_0(Y) \longrightarrow \tau_0(Z) \longrightarrow 0 \ . $$
These statements follow immediately from the long exact $Hom$-sequences
and the assumptions 2. and 3. Using the argument
of [KW], thm. 4.4
this implies

\begin{Lemma}. {\it Put $H^i(X)=\tau_0(X[i])$. Then for a distinguished
triangle $(X,Y,Z)$ in $\cal B$ there exists a long exact cohomology
sequence in ${\cal B}_0$}
$$ ... \to H^{-1}(Z) \to H^0(X) \to H^0(Y) \to H^0(Z) \to H^1(X) \to ... \ .$$
\end{Lemma}

\bigskip\noindent
In our case ${\cal B}={\cal B}_{m\vert n}$ satisfies these properties 1.-4., and indeed the suspension functor is induced by  
the parity shift functor $\Pi$ on $svec_K$ via the equivalence ${\cal B} \sim Rep_k(H) \otimes_k svec_K$ of categories.  Using
$A[2] \cong A$ defined by $id_A \otimes \zeta_2^{-1}$ and the identification ${\cal B} = sRep_K(H)$, we may identify the suspension 
functor with the parity shift functor, with functorial isomorphisms $H^{2i}(A)\cong A_+$ and $H^{2i+1}(A) \cong A_-$ for $A=A_+ \oplus A_-$ in $sRep_K(H)$.  Hence the long exact sequence
of the cohomology becomes a hexagon. ${\cal B}_{m\vert n}$ is 
a $\Pi$-category, a triangulated category enhanced by a super space structure, in the following sense: 

\bigskip\noindent
By definition a $\Pi$-category is a triangulated category with the properties 
1.-4. from above such that there exist functorial isomorphisms $\Pi^2(A)\cong A$ for
the suspension $\Pi(A) := A[1]$. For a functor $F:{\cal A} \to \cal B$ then notice
$F(A) \cong F(A)_+ \oplus F(A)_-$, where $H^0F(A)=F(A)_+$ and $H^1F(A)=F(A)_-$.
An additive functor $F$ from an abelian category $\cal A$  to a $\Pi$-category $\cal B$ will be called weakly exact, if $F(A)_+$ and $F(A)_-$ transform short exact sequences $0\to A \to B\to C \to 0$
into exact hexogons in ${\cal B}_0$ $$  \xymatrix{ 
 & F(B)_+ \ar[r]& F(C)_+ \ar[dr]&  \cr
 F(A)_+ \ar[ur]&  &  &   F(A)_- \ar[dl]\cr
& F(C)_- \ar[ul]& F(B)_- \ar[l]&   \cr}
$$
From the definition of the functor $\varphi = \gamma\circ \beta \circ \alpha$
the following then is obvious

\bigskip\noindent
\begin{Lemma} The functor $\varphi: {\cal T}_{m\vert n} \to sRep_K(H)$
is weakly exact. 
\end{Lemma}

\newpage
\centerline{\bf References}

\bigskip\noindent
[AK] Andre Y., Kahn B., Nilpotence, Radiceaux et structures monoidales, \hfill\break
arXiv/math/020327343 (2002)

\bigskip\noindent
[BR] Berele A., Regev A., Hook Young diagrams with applications to combinatorics and to representation theory of Lie superalgebras, Adv. in Math. 64 (1987), 118 - 175

\bigskip\noindent
[BKN1], Boe B.D., Kujawa J.R., Nakano D.K, Cohomology and support varieties for Lie superalgebras (2008)

\bigskip\noindent
[BKN2], Boe B.D., Kujawa J.R., Nakano D.K, Cohomology and support varieties for Lie superalgebras II (2008)

\bigskip\noindent
[Ba] Balmer P., Spectra, spectra, spectra - Tensor triangular spectra
versus Zariski spectra of endomorphism rings, preprint

%\bigskip\noindent
%[B] Batchelor M., A decomposition theorem for comodules, Compositio Math t.34, n 2 (1977)

\bigskip\noindent
[Br] Brundan J., Kazhdan-Lusztig polynomials and character
formulae for the Lie superalgebra $gl(m\vert n)$, Journal of the AMS, vol. 16, n.1, p. 185 -- 231
(2002)

\bigskip\noindent
[BS] Balmer P., Schlichting M., Idempotent completion of triangulated categories, J. Algebra, 236(2) (2001), 819 - 834

\bigskip\noindent
[BS1] Brundan J., Stroppel C., Highest weight categories arising from Khovanov's diagram
algebra I: Cellularity, arXiv (2009)

\bigskip\noindent
[BS2] Brundan J., Stroppel C., Highest weight categories arising from Khovanov's diagram
algebra II: Kostantity, arXiv (2009)

\bigskip\noindent
[BS3] Brundan J., Stroppel C., Highest weight categories arising from Khovanov's diagram
algebra III: Category $\cal O$, arXiv (2010)

\bigskip\noindent
[BS4] Brundan J., Stroppel C., Highest weight categories arising from Khovanov's diagram
algebra IV: The general linear supergroup, arXiv (2010)

%[BW] Brzezinski-Wisbauer, Corings and Comodules,
%London Math. Society Lecture Notes Series 309

%\bigskip\noindent
%[C] Chen, The stable monomorphism category of a Frobenius category

\bigskip\noindent
[D] Deligne P., Categories tensorielles, Moscow Mathematical Journal, vol 2, number 2 (2002), 227 - 248

\bigskip\noindent
[D2] Deligne P., Categories Tannakiennes, The Grothendieck Festschrift, vol II, Progress in Mathematics 87, Birkh\"auser 1990

\bigskip\noindent
[DM] Deligne P., Milne J.S., Tannakian Categories, in Hodge Cycles, Motives, and Shimura Varieties, Springer Lecture Notes 900, (1982)

\bigskip\noindent
[DS] Duflo M., Serganowa V., On associated variety for Lie Superalgebras, 2008

\bigskip\noindent
[DS] Dwyer W.G., Spalinski J., Homotopy theories and model categories,
Handbook of algebraic topology ed. I.M.James, Elsevier (1995)

%\bigskip\noindent
%[Fa] Faith C., Algebra II, Ring theory, Grundlehren der mathematischen Wissenschaften 191,
%Springer (1976)

\bigskip\noindent
[G] Garmonie J., Indecomposable representations of special linear Lie superalgebras,
J. of Algebra 209, 367 - 401 (1998)

\bigskip\noindent
[GM] Gelfand S.I., Manin.Y.I, Methods of Homological Algebra, Springer Monographs in Mathematics

\bigskip\noindent
[He1] Heidersdorf T., Semisimple quotients of representation categories
of Lie superalgebras and the case of sl(2,1),
preprint (2010)

\bigskip\noindent
[He2] Heidersdorf T., Representations of the Lie superalgebra osp(2,2n)
and the Lie algebra sp(2n-2),
preprint (2010)

\bigskip\noindent
[Ha] Happel D., Triangulated categories in the representation theory of finite dimensional
algebras, London Math Society Lect. series 119, Cambridge universty press (1988)

\bigskip\noindent
[Hi] Hirschhorn P.S., Model categories and their localizations, (2003), ISBN 0-8218-3279-4

\bigskip\noindent
[H] Hovey M., Palmieri J. H., Strickland N.P., Axomatic stable homotopy theory, Amer. Math. Soc. 128 (1997), no. 610, x + 114, AMS

\bigskip\noindent
[H] Hovey M., Model categories, Mathematical Surveys and Monographs vol 63, AMS

%\bigskip\noindent
%[K] Kan D.,

%\bigskip\noindent
%[HS] Hilton-Stammbach, A course in Homological Algebra, Graduate Texts in Mathematics, Springer
%1970

%\bigskip\noindent
%[H\"u] H\"uffmann Andreas, On representations of super coalgebras,
%arXiv:hep-th/9403100v1

\bigskip\noindent
[JHKTM] Van der Jeugt J., Hughes J.W.B., King R.C., Thierry-Mieg J., Character formulas
for irreducible modules of the Lie superalgebra $sl(m\vert n)$, J. Math. Phys. 31 (1990), no.9.
2278 -2304

\bigskip\noindent
[JHKTM2] Van der Jeugt J., Hughes J.W.B., King R.C., Thierry-Mieg J., A character formulas
for simply atypical modules of the Lie superalgebra $sl(m\vert n)$, Comm. Algebra 18 (1990), no.10, 3453 - 3480,
2278 -2304

\bigskip\noindent
[Ke] Keller B., Derived categories and tilting, in: Handbook of tilting theory, 49 - 104,
London Math Society Lect. series 332, Cambridge universty press (2007)

\bigskip\noindent
[Ke2] Keller B., Derived categories and their uses, preprint

\bigskip\noindent
[Li] Littlewood  D.E., The theory of group characters, Oxford university press, Oxford 1950

\bigskip\noindent
[Na] Nakaoka H., General heart construction on a triangulated category (I):
Unifying t-structures and cluster tilting subcategories, arXiv:0907.2080v6

\bigskip\noindent
[N] Neeman A., Triangulated categories, Annals of mathematics studies 148, Princeton university press 2001

%\bigskip\noindent
%[QS] Quella Th.-Schomerus V., Free fermion resolution of supergroup WZNW models,
%arXiv/0706.0744v1

%\bigskip\noindent
%[R] Roig A., Minimal resolutions and other minimal models, publicacions mathematiques, vol 37 (1993), 285 - 303

\bigskip\noindent
[Sch] Scheunert M., The Theory of Lie Superalgebras, Lecture Notes in Mathematics 716,
Springer 1979

\bigskip\noindent
[Se] Sergeev A.N., Tensor algebra of the identity representation
as a module over the Lie superalgebras $Gl(n,m)$ and $Q(n)$, Mat.
Sb 123 (165) (1984), no.3, 422 - 430

%\bigskip\noindent
%[S] Sweedler M.E., Hopf algebras, Benjamin 1969

%\bigskip\noindent
%[Si] Simson, On coalgebras of tame comodule type, in: Representations of algebras II. Bejing
%Normal University (2002)

%\bigskip\noindent
%[Su] Su Y., Composition Factors of Kac-modules for the general linear
%Lie superalgebras $gl(m\vert n)$

%\bigskip\noindent
%[SZ] Su Y.-Zhang R.B., Generalized Jantzen filtration of Lie superalgebras I

%

%\bigskip\noindent
%[T] Takeuchi, Morita theorems for categories of comodules, J. Fac. Sci. Univ. Tokyo 24 (1977), 629 - 644

\bigskip\noindent
[V] Verdier J.L., Categories derivees, in SGA 4$\frac{1}{2}$,
Springer Lecture Notes 569, Cohomologie Etale

\bigskip\noindent
[W] Weissauer R., Model structures, categorial quotients and representations of super commutative Hopf algebras I, preprint 2010

\bigskip\noindent
[W1] Weissauer R., Semisimple algebraic tensor categories,
arXiv:0909.1793v2

%\bigskip\noindent
%[WK] Kr\"amer Th.-Weissauer R.,  On the tensor square of irreducible representations
%of reductive Lie superalgebras, ArXiv:0910.5212v1

%\bigskip\noindent
%[We] Weselmann U., private communication

%\bigskip\noindent
%We remark in passing that similarly one can show that the functors $F'=U$
%$$ F'= Res_{p^+}^g: \cal C \to \cal D \ .$$
%$$  U'= CoInd_{p^+}^g : \cal D \to \cal C  $$
%define an adjoint pair of Quillen functors for the stable model structure
%on $\cal D$.
%

\end{document}